\documentclass{article}
\usepackage[utf8]{inputenc}
\usepackage{amsfonts}
\usepackage{amsmath}
\usepackage{amsthm}
\usepackage[colorlinks,allcolors=blue]{hyperref}
\usepackage[percent]{overpic}

\usepackage{algorithm,algpseudocode}

\usepackage[T1]{fontenc}

\usepackage{xcolor}

\newcommand{\figref}[1]{Figure~\ref{#1}}

\newcommand{\diag}[1]{\,\mathrm{diag}\left(#1\right)}

\newcommand{\R}{\mathbb{R}}
\newcommand{\C}{\mathbb{C}}
\def\ee{\ensuremath{\mathrm{e}}}
\def\ii{\ensuremath{\mathrm{i}}}
\DeclareMathOperator{\real}{Re}
\DeclareMathOperator{\imag}{Im}
\DeclareMathOperator{\conj}{conj}

\DeclareMathOperator*{\argmax}{arg\,max}
\DeclareMathOperator*{\argmin}{arg\,min}

\newtheorem{theorem}{Theorem}
\newtheorem{remark}[theorem]{Remark}
\newtheorem{corollary}[theorem]{Corollary}
\newtheorem{proposition}[theorem]{Proposition}

\numberwithin{equation}{section}

\title{Unitarity of some barycentric rational approximants}
\author{Tobias Jawecki\footnote{
Institute for Theoretical Physics, Vienna University of Technology,
Vienna, Austria,
\href{mailto:tobias.jawecki@tuwien.ac.at}{tobias.jawecki@tuwien.ac.at}
}~~and
Pranav Singh\footnote{
Department of Mathematical Sciences,
University of Bath,
Bath, United Kingdom,
\href{mailto:ps2106@bath.ac.uk}{ps2106@bath.ac.uk}
}}
\date{May 2022}

\begin{document}

\maketitle
\begin{abstract}
The exponential function maps the imaginary axis to the unit circle and, for many applications,
this unitarity property is also desirable from its approximations.
We show that this property is conserved not only by the $(k,k)$-rational barycentric interpolant of the exponential on the imaginary axis,
but also by $(k,k)$-rational barycentric approximants that minimize a linearized approximation error.
These results are a consequence of certain properties of singular vectors of Loewner-type matrices associated to linearized approximation errors.
Prominent representatives of this class are rational approximants computed
by the adaptive Antoulas--Anderson (AAA) method and the AAA--Lawson method.
Our results also lead to a modified procedure with improved numerical stability of the unitarity property and reduced computational cost.
\end{abstract}

% 15A23 Factorization of matrices
% 41A20 Approximation by rational functions
% 65D15 Algorithms for approximation of functions
% these subject classifications are also similar to what is used in AAA and AAA-Lawson
% ACM
% G.1.1: Interpolation
% G.1.2: Approximation
\section{Introduction}

Polynomial and rational approximations to functions have a wide range of applications \cite{TrefethenATAP}.
Rational approximations have some advantages compared to polynomial approximation \cite{salzer1981rational,schneider1986some, berrut1988rational},
two main strengths being
strong performance for approximating functions near singularities
and on unbounded domains~\cite{TrefethenATAP}.
%(i) for the approximation of functions with singularities, a rational approximant with similar
%singularities can perform well \cite{baker_graves-morris_1996}, and (ii) rational approximants with singularities outside of
%the domain of interest are bounded in norm.
% this holds for $(s,k)$-rational function with $s\leq k$.

In the present work we are concerned with $(k,k)$-rational approximants to the imaginary exponential,
\begin{equation}
    \label{eq:ratapprox}
    r(x) = \frac{p(x)}{q(x)} \approx \ee^{\ii x}, \qquad x \in I \subset \mathbb{R},
\end{equation}
where $I$ is a bounded interval in $\mathbb{R}$, and $p(x)$ and $q(x)$
are complex polynomials of degree $k$ each.
This is closely related to the time integration of ordinary differential equations (ODEs), and different approximations to the exponential function yield different numerical methods for time integration. In this context, boundedness of the underlying rational function on a specific subset of the complex plane results in stable numerical integrators \cite{HW02}.

In particular, in context of the application specified in~\eqref{eq:ratapprox}, i.e. the approximation of the exponential function on an interval on the imaginary axis, the main advantage of rational approximation over polynomial approximation is that there exist rational approximants that satisfy the unitarity property,
\begin{equation}
    \tag{$\star$} \label{eq:unitarity}
|r(x)| = 1, \qquad x \in \mathbb{R}.
\end{equation}
Unitarity of the rational approximants to the imaginary exponential
has strong benefits for ODEs with a skew-Hermitian structure \cite{HWL02}. The requirement for unitarity~\eqref{eq:unitarity} often arises in the context of equations of quantum mechanics \cite{lubich08fqc, faou12gni}, for instance.

An important class of $(k,k)$-rational approximants which satisfy the unitarity property~\eqref{eq:unitarity} are (diagonal) Pad\'{e} approximants \cite{baker_graves-morris_1996}.
Pad\'{e} approximants have a high order of accuracy around the origin. However, being Taylor based and thus asymptotic in nature, they prove inadequate
when a more uniform accuracy over the interval $I$ in~\eqref{eq:ratapprox} is desired, or other interpolation properties have to be satisfied.

\begin{figure}
\centering
\begin{overpic}
[width=0.8\textwidth]{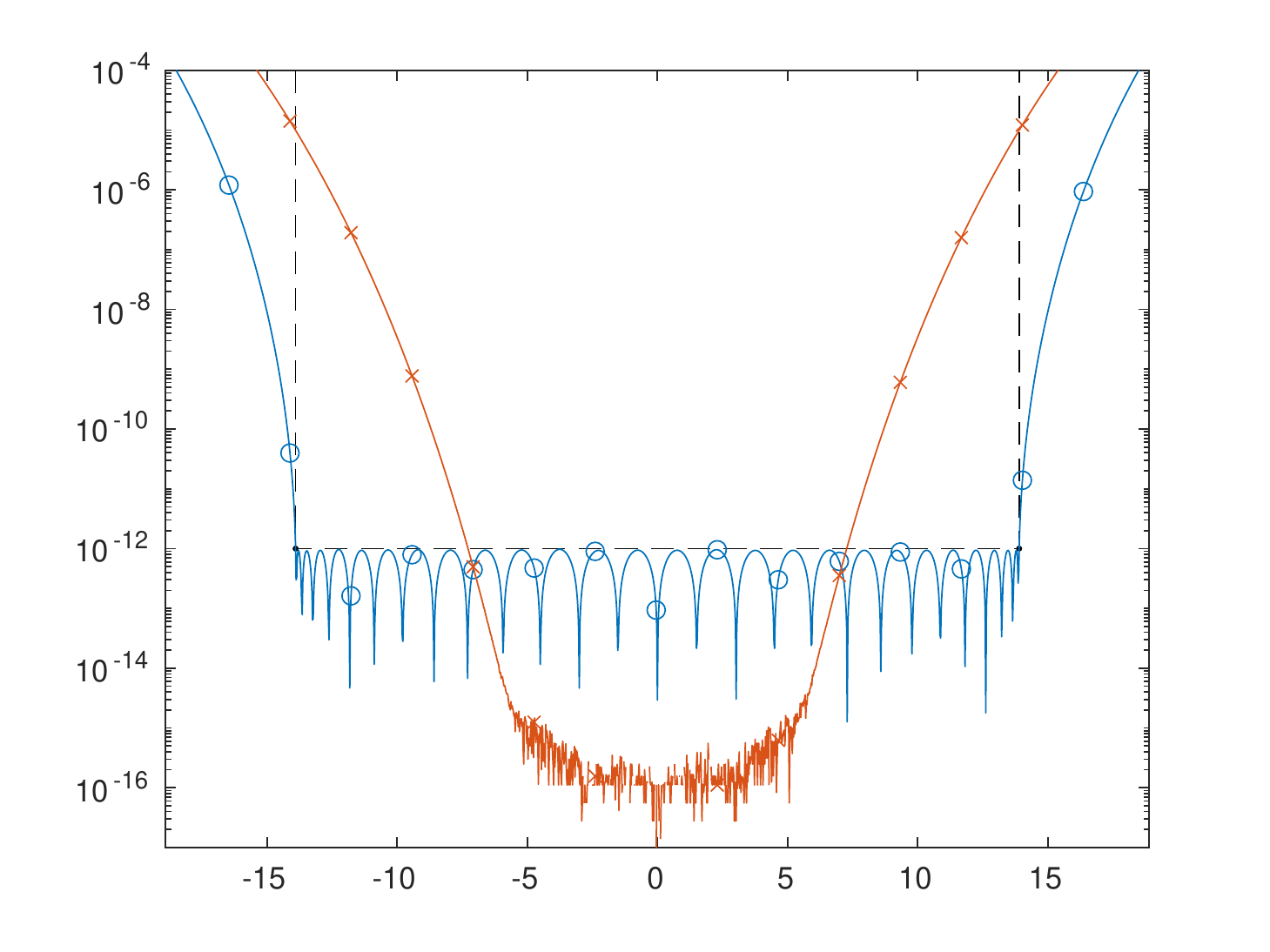}
\put(-11,40){$|r(x)-\ee^{\ii x}|$}
\put(50,0){$x$}
\end{overpic}
\caption{This figure shows the error in two different rational approximations
to the imaginary exponential. The error of the diagonal Pad\'{e} approximant of degree~$13$
is marked by the symbol~($\times$), while the symbols~($\circ$) mark the error of a rational approximant
generated by the AAA--Lawson method.
Namely, the AAA--Lawson method is applied to approximate
the imaginary exponential on the interval~$[-13.9,13.9]$
using a~$(13,13)$-rational function.
The error of the AAA--Lawson method is less than~$10^{-12}$
uniformly on the interval, and this error bound is illustrated by a dashed line.
}
\label{fig.AAAnPade}
\end{figure}

A more flexible approximation is provided by rational fitting algorithms such as the adaptive Antoulas--Anderson (AAA) method \cite{NST18} and the AAA--Lawson method~\cite{NT20}, which can provide more accuracy over a relevant interval $I$ (or a discrete set of points) than a diagonal Pad\'{e} approximant of the same degree. For instance, in \figref{fig.AAAnPade} we find that a $(13,13)$-rational approximant produced by the AAA--Lawson algorithm provides uniformly high accuracy ($<10^{-12}$) over the interval of interest.  While the diagonal Pad\'{e} approximant of degree $13$ has a very high accuracy in the neighborhood of $x=0$, its accuracy towards the extremes of the interval is substantially lower (roughly $10^{-5}$ near $x=-13.9$, for instance). % AAA similar

In the present paper, we consider $(k,k)$-rational approximants in barycentric rational representation
of three types:
\begin{enumerate}
    \item[(i)] rational interpolants of $\ee^{\ii x}$ at exactly $2k+1$ nodes,
    \item[(ii)] rational approximants of $\ee^{\ii x}$ that interpolate $\ee^{\ii x}$ at $k+1$
{\em support nodes} and minimize a linearized error at a larger number of so called {\em test nodes}, and
    \item[(iii)] non-interpolatory rational approximants of $\ee^{\ii x}$ that minimize a linearized error at test nodes.
\end{enumerate}
In particular, rational approximants produced by the AAA algorithms fall under type (ii). AAA can also generate approximants with uniform accuracy over an interval $I$ by adaptively choosing (test and support) nodes.
Rational approximants generated by the AAA--Lawson algorithm fall under type (iii).
The AAA--Lawson algorithm further increases the accuracy of the approximant generated by AAA
by iteratively solving re-weighted least-squares problems.
Using a barycentric rational representation for the rational approximation yields two benefits. The first advantage of this representation is
strong stability properties~\cite[Subsection~2.3]{FNTB18}.
The other advantage is that the rational approximants of types (i)-(iii) can be found by computing a singular
value decomposition of a Loewner matrix \cite{Be70,berrut2000matrix}
or an expanded Loewner matrix~\cite{NT20}.

We show that, while being much more flexible, rational approximants of the type (i)--(iii) share an important property with diagonal Pad\'{e} approximants: Unitarity~\eqref{eq:unitarity}, which is highly desirable for applications.
In contrast, while the \texttt{rkfit} procedure~\cite{BG15} can produce rational approximants with an accuracy comparable to the AAA and AAA--Lawson methods, the approximants generated by \texttt{rkfit} are not unitary in general. This makes AAA and AAA--Lawson methods better suited to applications in~\eqref{eq:ratapprox}.

\paragraph{Outline of the paper.}

% OVERVIEW ON LABELS
% 2 Section~\ref{sec:barfct}
% Proposition~\ref{prop.2mm1gencaesunitary}
% 3 Section~\ref{sec:LoewnerMat}
% 3.1 Subsection~\ref{subsec:caseI}
% Corollary~\ref{prop.main1}
% 3.2 Subsection~\ref{subsec:caseII}
% 3.3 Subsection~\ref{subsec:rescaleLoewnerandunitary}
% Proposition~\ref{prop.rminisuunitary}
% 3.4 Subsection~\ref{subsec:caseIImu}
% 4 Section~\ref{sec:caseIII}
% Proposition~\ref{prop:maincaseIII}
% 5 Section~\ref{subsection:Cayleyform}
% Algorithm~\ref{alg:modAAA}
% Algorithm~\ref{alg:modAAALawson}
% 6 Section~\ref{sec:AAA}
% A Appendix~\ref{subsec.Asvd}
% A.1 Appendix~\ref{subsec.appendixcaseIII}
% B Appendix~\ref{subsec.Aauxiliaryunitary}
In Section~\ref{sec:barfct} we recall barycentric rational representations and show, in Proposition~\ref{prop.2mm1gencaesunitary}, that rational interpolants that interpolate $\ee^{\ii x}$ at a maximal number of distinct nodes -- i.e., of type (i) above -- are unitary. 

In Section~\ref{sec:LoewnerMat}, we consider an interpolatory barycentric representation, and express a linearized error in terms of a Loewner matrix.  
In Subsection~\ref{subsec:caseI}, we consider the case where the linearized error is zero.
This corresponds to barycentric rational interpolation which falls under type~(i), and consequently, unitarity follows from Proposition~\ref{prop.2mm1gencaesunitary}.
% which falls under type (i) and consequently unitarity follows from  Proposition~\ref{prop.2mm1gencaesunitary}. 

The remaining manuscript is concerned with types (ii) and (iii) -- i.e., rational approximants that minimize a linearized error. In Subsection~\ref{subsec:caseII}, we consider interpolatory barycentric approximants of type (ii).
Specifically, the approximants generated by the AAA method fall in this class. 
In Proposition~\ref{prop.rminisuunitary}, Subsection~\ref{subsec:rescaleLoewnerandunitary},
we show that approximants in this class are unitary. These results are generalized to the case of a weighted linearized error in Subsection~\ref{subsec:caseIImu}.
In Section~\ref{sec:caseIII} we consider approximants in a
non-interpolatory barycentric representation -- i.e., type (iii) -- as used in the AAA--Lawson method.
Unitarity for such approximants is shown in Proposition~\ref{prop:maincaseIII}.

\begin{figure}
\centering
\begin{overpic}
[width=0.8\textwidth]{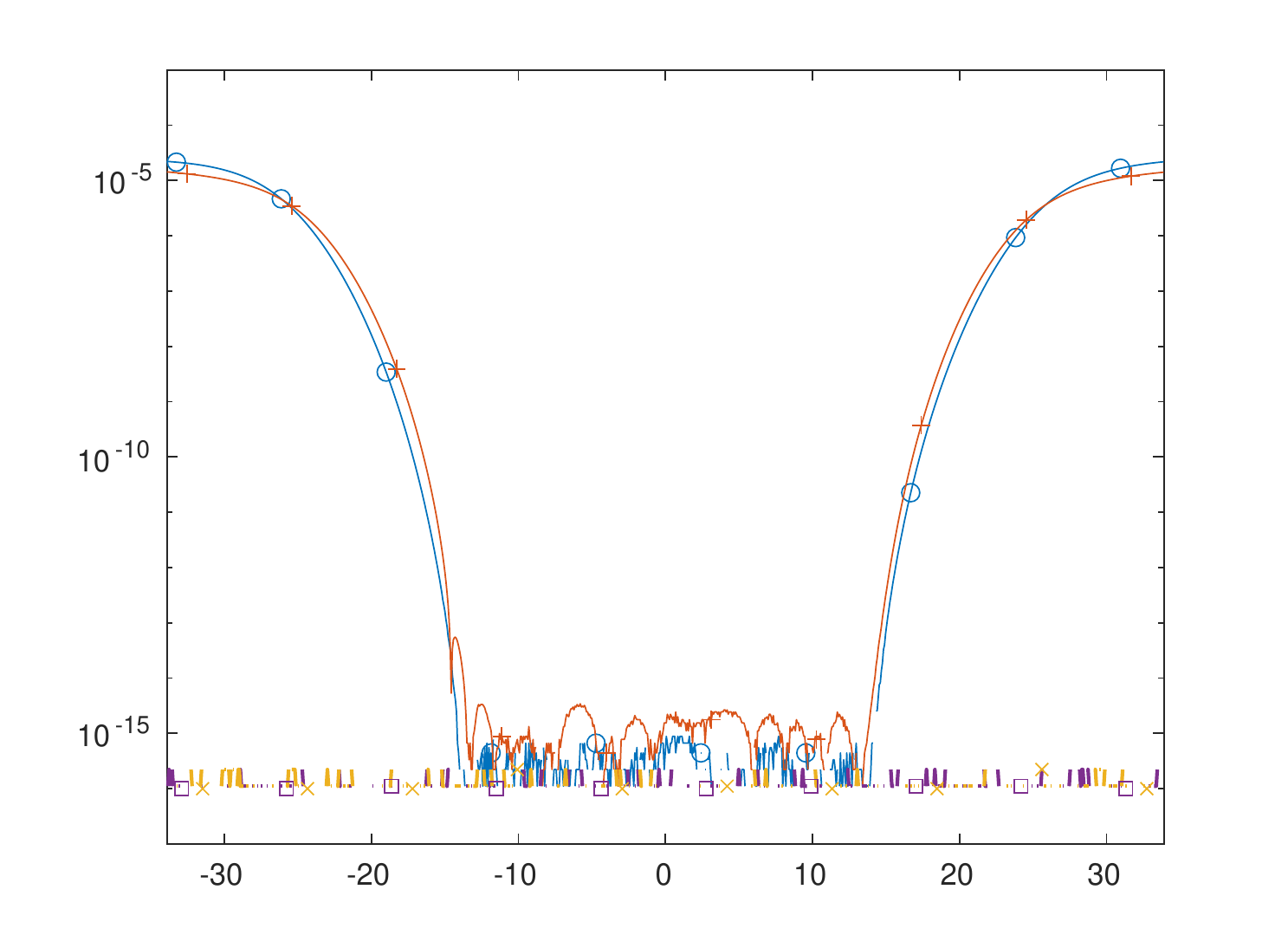}
\put(-11,40){$||r(x)|-1|$}
\put(50,0){$x$}
\end{overpic}
\caption{This figure shows the deviation of the
unitarity property in computer arithmetic for
different barycentric rational approximants.
Namely, for the
barycentric rational approximants $r$ generated by the original~($\circ$) and modified~($\square$) AAA algorithms, and $r_b$ generated by the original~($+$) and modified~($\times$) AAA--Lawson algorithms.
The modified algorithms are introduced in Section~\ref{subsection:Cayleyform}.
The approximants $r$ and $r_b$ are $(14,14)$ and $(13,13)$-rational functions, respectively,
and all approximants are generated to approximate $\ee^{\ii x}$
for $x\in[-13.9,13.9]$ with an approximation
error $\leq 10^{-12}$ on this interval.
The approximation error of $r_b$ is also illustrated in~\figref{fig.AAAnPade}.
All of these approximants are unitary in theory,
however, in computer arithmetic this property is not exactly preserved
when using the original AAA and AAA--Lawson methods. In particular, the deviation in unitarity is below $10^{-14}$ on the domain of approximation, $x\in[-13.9,13.9]$, but at $x=35$ it is nearly %$10^{-2}$ and $10^{-4}$, respectively, 
$10^{-5}$ for the original AAA and AAA-Lawson methods.
Unitarity can be maintained uniformly on $\R$ in computer arithmetic by using the
modified algorithms -- particularly, the deviation of the unitarity property for the modified algorithms is approximately $10^{-16}$ or exactly zero (not visible in this figure).}
\label{fig.UnitarityRoundoff}
\end{figure}

Despite the theoretical unitarity of rational approximants of types (ii) and (iii), due to the finite precision of computer arithmetic, in practice these approximants tend to deviate significantly from unitarity away from the domain of approximation (see  \figref{fig.UnitarityRoundoff}). 
To remedy this situation, we describe a slight modification to the original AAA and AAA--Lawson algorithms in Section~\ref{subsection:Cayleyform}. In particular, our approach replaces complex SVD with real SVD, reducing computational cost, and we resort to a Cayley representation, expressing the numerator as the complex conjugate of the denominator. In particular, as illustrated in \figref{fig.UnitarityRoundoff}, the modified algorithms,  Algorithm~\ref{alg:modAAA} and~\ref{alg:modAAALawson}, show unitarity at machine precision even away from the domain of approximation. 

In Section~\ref{sec:AAA} we briefly sketch the AAA and AAA--Lawson method
to illustrate that unitarity indeed holds for approximants generated
by these methods. The unitarity results presented in the present paper are based on certain properties of singular vectors
of Loewner-type matrices, which are derived in Appendix~\ref{subsec.Asvd}.

\section{Barycentric rational representation}
\label{sec:barfct}

In the present work we consider rational approximations
to the imaginary exponential function.
In particular, we make use of barycentric rational representations of rational functions.
Approximants in the representations~\eqref{eq.defr} and~\eqref{eq.defr2} below are topic
of Section~\ref{sec:LoewnerMat} and~\ref{sec:caseIII}, respectively.

The following barycentric rational representation relies on distinct {\em support nodes} $y_1,\ldots,y_m\in\R$
and coefficients $w_1,\ldots,w_m\in\C^m$.
For the imaginary exponential function evaluated at the support nodes we also use the notation $f_j=\ee^{\ii y_j}$.
In Section~\ref{sec:LoewnerMat} we consider barycentric rational representations
given by the quotient of partial fractions
\begin{equation}\label{eq.defr}
r(x) = \sum_{j=1}^m \frac{f_j \,w_j}{x-y_j}
\Big/ \sum_{j=1}^m \frac{w_j}{x-y_j} =: n(x)/d(x).
\end{equation}
The poles of the partial fractions $n$ and $d$ coincide with the support nodes.
However, for the function $r$ these singularities are removable.
Assuming $w_j\neq 0$, we have the identity
\begin{equation}\label{eq.limrzj}
r(y_j) = \ee^{\ii y_j},~~~j=1,\ldots,m,
\end{equation}
as a limit, so that $r$ interpolates the imaginary exponential at the support nodes.

In the sequel, the notation $(m-1,m-1)$-rational function refers
to a rational function $r=p/q$,
where $p$ and $q$ are polynomials of degree $\leq m-1$.
A function in barycentric rational representation~\eqref{eq.defr}
with $m$ support nodes
corresponds to an $(m-1,m-1)$-rational function. Let
$$
\ell(x) = \prod_{j=1}^m (x-y_j)
$$
denote the polynomial with zeros located at the poles of
the partial fractions $n$ and $d$.
We introduce functions~$p$ and~$q$ as
\begin{equation}\label{eq.defpq}
\begin{aligned}
&p(x) := \ell(x)n(x)
= \sum_{j=1}^m f_j w_j \prod_{k\neq j} (x-y_k),
~~~\text{and}\\
&q(x) := \ell(x)d(x)
= \sum_{j=1}^m w_j \prod_{k\neq j} (x-y_k).
\end{aligned}
\end{equation}
Indeed, $p$ and $q$ correspond to polynomials of degree $\leq m-1$.
Multiplying the numerator and denominator of the barycentric rational representation $r$ by $\ell$,
we observe
\begin{equation}\label{eq.aaarpq}
r(x) = \frac{\ell(x)\, n(x)}{\ell(x)\, d(x)} = \frac{p(x)}{q(x)},
\end{equation}
which shows that $r$ is an $(m-1,m-1)$-rational function.

For non-interpolatory approximants in Section~\ref{sec:caseIII}
we also consider barycentric rational representations
\begin{equation}\label{eq.defr2}
r_b(x) = \sum_{j=1}^m \frac{\alpha_j}{x-y_j}
\Big/ \sum_{j=1}^m \frac{\beta_j}{x-y_j},
\end{equation}
where $\alpha_1,\ldots,\alpha_m\in\C$ and
$\beta_1,\ldots,\beta_m\in\C$ denote coefficients,
and $y_1,\ldots,y_m\in\R$ denote support nodes.
In the case of $\alpha_j=f_j w_j$ and $\beta_j=w_j$, this representation
coincides with~\eqref{eq.defr}.
In contrast to~$r$ in~\eqref{eq.defr}, $r_b$ does not necessarily interpolate
$\ee^{\ii x}$ at the support nodes~\eqref{eq.limrzj}.
In particular, this interpolation property does not hold
if $\alpha_j\neq f_j\beta_j $.
Expanding $r_b$ similar to~\eqref{eq.aaarpq},
shows that $r_b$
corresponds to an $(m-1,m-1)$-rational function.

We refer to a rational function $r=p/q$ as \textit{irreducible}
if the polynomials $p$ and $q$ have no common zeros.
Furthermore, the \textit{poles} of $r$ refer to the zeros of~$q$.
These terms are used in an equivalent manner for a rational function
given in barycentric rational representation.
%based on the representation $r=p/q$ as in~\eqref{eq.aaarpq}.

We proceed with some general results concerning $(m-1,m-1)$-rational interpolations to the imaginary exponential function.
Similar results are specified for
barycentric rational interpolation in Subsection~\ref{subsec:caseII}.

\paragraph{General remarks on interpolation
by~\texorpdfstring{$(m-1,m-1)$}{(m-1,m-1)}-rational functions.}
\label{subsec.generalcase}

For an $(m-1,m-1)$-rational function $r=p/q$,
the numerator and denominator polynomials are of degree $\leq m-1$ and
have $m$ coefficients each.
However, the resulting rational function only has $2m-1$ free parameters
(this result can be observed via partial fraction decomposition).
It is natural  to consider $(m-1,m-1)$-rational functions
to interpolate the imaginary exponential at $2m-1$ given nodes
in a general setting.
At this point we assume that such a rational interpolation
function exists and is irreducible,
for further details we refer to~\cite[Section 2]{Be70}.
The following proposition shows unitarity for this class of rational interpolants.

\begin{proposition}\label{prop.2mm1gencaesunitary}
Let $r$ be an irreducible $(m-1,m-1)$-rational function,
and let $\theta_1,\ldots, \theta_{2m-1}\in\R$
be distinct nodes such that $r$ satisfies the interpolation property
\begin{equation}\label{eq.rinterpolexp2mm1}
r(\theta_j)=\ee^{\ii \theta_j},~~~j=1,\ldots,2m-1.
\end{equation}
Then $r$ has no poles on the real axis, and satisfies
\begin{equation}
\tag{\ref{eq:unitarity}} \label{eq.rinterpolexp2mm1x}
|r(x)|=1,~~~\text{for}~~x\in\R.
\end{equation}
\end{proposition}
A proof of this proposition is provided in Appendix~\ref{subsec.Aauxiliaryunitary}.

%%%%%%%%%%%%%%%%
%%%%%%%%%%%%%%%%

%%%%%%%%%%%%%%%%
%%%%%%%%%%%%%%%%

\section{Interpolatory rational approximation and Loew\-ner matrices}
\label{sec:LoewnerMat}

We proceed to consider rational approximants of $\ee^{\ii x}$ in
barycentric rational representation~\eqref{eq.defr}.
In the present section, we consider $(m-1,m-1)$-rational approximants where $m$ is fixed.
Approximants of a moderate degree are favourable due to various reasons,
e.g., computational cost.
Furthermore, a certain approximation accuracy over a given interval or discrete set of nodes
is desirable in practice.
Due to limitations of the methods discussed in the present work,
our focus is on the discrete case.
This includes discretized intervals,~e.g.,~\cite{NST18,NT20}.
%The methods discussed in the present work cannot provide accuracy on an interval in a direct manner.
%Instead, we assume intervals to be discretized by a sufficiently large number of nodes.
% set of nodes can be chosen in an adaptive manner~\cite{NST18,NT20}.
We introduce the notation $x_1,\ldots,x_n\in\R$ for the {\em test nodes},
over which the rational approximant $r$ needs to approximate $\ee^{\ii x}$,
\begin{equation}\label{eq:approxxk}
r(x_k)\approx \ee^{\ii x_k},~~~k=1,\ldots,n,
\end{equation}
with a prescribed accuracy.
In addition, the interpolation property at the support nodes~\eqref{eq.limrzj}
is desirable in an interpolatory setting.
In the sequel, we assume that the support nodes~$y_1,\ldots,y_m\in\R$
and the test nodes~$x_1,\ldots,x_n\in\R$ are disjoint sets
of distinct nodes, unless explicitly stated otherwise.

A special case occurs for $n=m-1$,~i.e.,~we have $m$ given support
nodes and $m-1$ given test nodes.
This yields a total of $2m-1$ nodes,
and as remarked at the end of the previous section,
rational interpolation is viable in this setting.
Similar interpolation problems have been studied earlier in~\cite{Be70,AA86,berrut2000matrix}.
%we recall the required theory in the following.
As a new result, we prove unitarity of such interpolants
in Corollary~\ref{prop.main1} in Subsection~\ref{subsec:caseI}
based on Proposition~\ref{prop.2mm1gencaesunitary}.
% which has not been stated before in this setting.
%Interpolation at preassigned test and support is also suggested

If a larger number of test nodes is given (namely, $n>m-1$),
then interpolation at all test and support nodes is
no longer viable in general.
In this case, we replace the interpolation property
at the test nodes by a near-best approximation property.
In particular, we consider barycentric rational approximants
that minimize a linearized error
over the test nodes, see Subsection~\ref{subsec:caseII}.
This approach is also utilized in the AAA method~\cite{NST18}.
For the case $n>m-1$, Proposition~\ref{prop.2mm1gencaesunitary}
does not apply due to non-interpolatory nature.
However, in Subsection~\ref{subsec:rescaleLoewnerandunitary}
we prove that the approximants that minimize the linearized error
are unitary -- a new result which is based on properties of singular vectors
of a re-scaled Loewner matrix.

We remark that the AAA method also includes an outer iteration to determine
the support nodes $y_1,\ldots,y_m$.
However, in this paper we show unitarity for arbitrary support nodes, and therefore we may consider them fixed. 
A sketch of the AAA method is also given
in Section~\ref{sec:AAA} further below.

Barycentric rational approximants that minimize a weighted linearized error
satisfy similar unitarity properties, as shown in Subsection~\ref{subsec:caseIImu}.
This setting is motivated by~\cite[Subsection~10]{NST18}, and its applications.
Weighted problems also appear in~\cite{BG15,BG17b}.

Non-interpolatory near-best approximants using the representation $r_b$ in~\eqref{eq.defr2} are discussed in Section~\ref{sec:caseIII} further below.
Approximants therein are related to the AAA--Lawson method~\cite{NT20}.

\medskip
We proceed to introduce common notation concerning the Loewner matrix,
which is relevant for barycentric rational approximation in the present section.

\paragraph{Loewner matrix.}
Assuming $y_1,\ldots,y_m\in\R$ are given support nodes
and $x_1,\ldots,x_n\in\R$ are given test nodes,
we define the Loewner matrix
\begin{equation}\label{eq.defLentries}
L\in\C^{n\times m},~~~\text{with} ~~ L_{kj} = \frac{\ee^{\ii x_k}-\ee^{\ii y_j}}{x_k-y_j}, 
\end{equation}
$k \in \{1,\ldots, n\},\ j \in \{1, \ldots, m\}$.
The Loewner matrix can be written as a product of matrices:
Let $C\in\R^{n\times m}$ denote the Cauchy matrix
\begin{subequations}\label{eq.defCauchyandDiag}
\begin{equation}\label{eq.defCauchyMat}
C\in\R^{n\times m},~~~\text{with}~~C_{kj} = \frac{1}{x_k-y_j},
\end{equation}
and let $S_F$ and $S_f$ denote the diagonal matrices
\begin{equation}\label{eq.defSFSf}
S_F=\diag{\ee^{\ii x_j}}\in\C^{n\times n},~~~
\text{and}~~~S_f=\diag{\ee^{\ii y_j}}\in\C^{m\times m}.
\end{equation}
\end{subequations}
Then the Loewner matrix can be expressed as
\begin{equation}\label{eq.LoewnerMatrixform}
L = S_F C - C S_f\ \ \in\C^{n \times m}.
\end{equation}
For the given support nodes, we consider a barycentric rational approximant
$r(x)=n(x)/d(x)$ as in~\eqref{eq.defr}.
We consider $w=(w_1,\ldots,w_m)^\top\in\C^m$ to be the vector
of the underlying coefficients in~\eqref{eq.defr}, which have not been
specified yet.
The matrices in~\eqref{eq.defCauchyandDiag} satisfy the following identities
when applied to $w$,
\begin{equation}\label{eq.SCapplyw}
\begin{aligned}
&(S_F Cw)_k = \ee^{\ii x_k} \sum_{j=1}^m \frac{w_j}{x_k-y_j} = \ee^{\ii x_k} d(x_k),
~~~\text{and},\\
&(C S_f w)_k = \sum_{j=1}^m \frac{f_j\, w_j}{x_k-y_j} = n(x_k),~~~\text{where}~~f_j=\ee^{\ii y_j}.
\end{aligned}
\end{equation}
Making use of the representation~\eqref{eq.LoewnerMatrixform}
and the identities~\eqref{eq.SCapplyw},
we conclude
\begin{equation}\label{eq.Lwtoxj}
(Lw)_k
= \ee^{\ii x_k}d(x_k) - n(x_k).
\end{equation}

%%%%%%%%%%%%%%%%%%%%%%%%%%%
\subsection{An \texorpdfstring{$(m-1,m-1)$}{(m-1,m-1)} rational interpolation
at \texorpdfstring{$2m-1$}{2m-1} nodes}
\label{subsec:caseI}

We proceed with an overview on barycentric rational interpolation
of $\ee^{\ii x}$ at given nodes.
Similar interpolation problems are studied earlier in~\cite{Be70,AA86,berrut2000matrix} and others.
We assume that $2m-1$ distinct nodes are given; specifically, we assume that we are given $m$ support nodes and $m-1$ test nodes,~i.e.,~we have the case $n=m-1$.
For an example concerning
interpolation of~$\ee^{\ii x}$ at preassigned test and support nodes
we also refer to~\cite[Subsection~III.B]{Kno08}.

In the present setting the Loewner matrix $L$ given in~\eqref{eq.defLentries}
has the dimension $m-1\times m$. Thus, $L$ has a non-trivial nullspace.
Let the vector of coefficients $w=(w_1,\ldots,w_m)^\top$
be in the nullspace of $L$,~i.e.,
$$
Lw=0.
$$
Then, due to~\eqref{eq.Lwtoxj}, the partial fractions
$n$ and $d$ satisfy
\begin{equation}\label{eq.Lwzero2}
n(x_k) = d(x_k) \ee^{\ii x_k},~~~k=1,\ldots,m-1.
\end{equation}
Furthermore, if the resulting barycentric rational approximation $r$
is irreducible, then~\eqref{eq.Lwzero2} implies
\begin{equation}
\tag{\ref{eq:unitarity}} \label{eq.m1mratxj}
r(x_k)=\ee^{\ii x_k},~~~k=1,\ldots,m-1.
\end{equation}

In the following corollary,
we summarize the interpolation properties of $r$,
together with a unitarity property which follows
Proposition~\ref{prop.2mm1gencaesunitary}.
\begin{corollary}
%[Corollary of Proposition~\ref{prop.2mm1gencaesunitary}]
\label{prop.main1}
Let $y_1,\ldots,y_m\in\R$ and $x_1,\ldots,x_{m-1}\in\R$ be given support and test nodes, respectively, $L$ be the corresponding Loewner matrix~\eqref{eq.defLentries},
and let $w=(w_1,\ldots,w_m)^\top\in\C^m$ denote a vector of coefficients that satisfies $Lw=0$.
Assume $w_j\neq 0$, and assume the generated barycentric rational function $r$
is irreducible.
Then $r$ interpolates $\ee^{\ii x}$
at the support and test nodes,
has no poles on the real axis, and satisfies
\begin{equation}
\tag{\ref{eq:unitarity}} |r(x)|=1,~~~x\in\R.
\end{equation}
\end{corollary}
\begin{proof}
The barycentric rational function $r$
interpolates at the support nodes due to~\eqref{eq.limrzj}
and the assumption $w_j\neq 0$.
Following~\eqref{eq.m1mratxj}, % and further remarks above,
the choice of $w$ implies that $r$ interpolates at the
test nodes.
Thus, $r$ interpolates at a total of $2m-1$ distinct nodes.
By construction, $r$ is an $(m-1,m-1)$-rational function.
Following Proposition~\ref{prop.2mm1gencaesunitary},
the interpolation properties of $r$
imply that $r$ has to no poles on the real axis and is unitary.
\end{proof}

In Subsection~\ref{subsec:rescaleLoewnerandunitary} further below,
we provide an alternative proof that the interpolant $r$ considered in Corollary~\ref{prop.main1} is unitary.
This involves a Cayley-type representation for $r$.

\subsection{Minimizing a linearized error}
\label{subsec:caseII}

We consider test nodes $x_1,\ldots,x_n\in\R$,
where, in contrast to the previous subsection, $n>m-1$.
%where $m$ is the number of given support nodes.
In the present subsection, we replace the interpolation property at the test nodes
by a near-best approximation property.
Namely, we aim to minimize a linearized error over the test nodes;
a practical approach which is utilized in the AAA method~\cite{NST18}.

%\begin{equation}\label{eq.approxexp}
%\ee^{\ii x} \approx r(x) = n(x)/d(x).
%\end{equation}
We recall $r(x) = n(x)/d(x) \approx \ee^{\ii x} $,
and we linearize this approximation property by multiplying by~$d(x)$ on both sides,
\begin{equation}\label{eq.approxlin}
n(x) \approx \ee^{\ii x} d(x).
\end{equation}

For the underlying test and support nodes we consider the Loewner matrix $L\in\C^{n\times m}$ as given in~\eqref{eq.defLentries}.
Let $w=(w_1,\ldots,w_m)^\top\in\C^m$ refer to the vector of coefficients of
$r$ as in~\eqref{eq.defr}.
Following~\eqref{eq.Lwtoxj}, the entries of $Lw$ evaluate the deviation in~\eqref{eq.approxlin} at the test nodes,
and the Euclidean norm of the vector $Lw$ yields
\begin{equation}\label{eq.linerrLWnorm}
\| Lw \|_2 = \left(\sum_{k=1}^n
|\ee^{\ii x_k} d(x_k) - n(x_k)|^2 \right)^{1/2}.
\end{equation}
This representation is referred to as the {\em linearized error}.

We aim to choose coefficients $w=(w_1,\ldots,w_m)^\top\in\C^m$
with $\|w\|_2=1$ such that the linearized error~\eqref{eq.linerrLWnorm} is minimized,~i.e.,
\begin{equation}\label{eq.minLw}
w=\argmin_{u\in\C^m,~\|u\|_2=1}\| Lu \|_2.
\end{equation}
Such coefficients are accessible by exploiting the
singular value decomposition of the Loewner matrix $L$.
We have the factorization
\begin{equation}\label{eq.svdL}
LV = US,
\end{equation}
where $V\in\C^{m \times m}$ and $U\in\C^{n \times m}$ refer to orthonormal bases of right and left singular vectors, respectively,
and $S=\diag{\sigma_1,\ldots,\sigma_m}\in\R^{m\times m}$ refers to
a diagonal matrix containing singular values of $L$.
We assume the ordering
$$
\sigma_1 \geq \ldots \geq \sigma_m \geq 0.
$$
A weight vector $w$ which minimizes the
linearized error norm as in~\eqref{eq.minLw}
is attained by a right singular vector
of $L$,~e.g.,\footnote{
  We use the notation $e_k=(0,\ldots,0,1)^\top\in\C^k$.
}~$w=V e_m$ with $\|Lw\|_2=\sigma_m$.
In particular, the smallest singular value of $L$ satisfies
$$
 \sigma_m \leq \|L u\|_2~~~\text{for any $u\in\C^m$ with $\|u\|_2=1$}.
$$

The generated barycentric rational approximation $r$ interpolates
at all test nodes if and only if $\sigma_m=0$.
%Thus, in that case $w$ as in~\eqref{eq.minLw} is in the nullspace
%of $L$, and $r$ is unitary following Corollary~\ref{prop.main1}
%in the previous subsection.
In general we have the case $\sigma_m>0$,
and interpolation at test nodes is not guaranteed.
Thus, Proposition~\ref{prop.2mm1gencaesunitary} does not apply.
Nevertheless, unitarity of the barycentric rational approximant
which minimizes the linearized error
is shown in the following subsection.

\subsection{A re-scaled Loewner matrix and unitarity}
\label{subsec:rescaleLoewnerandunitary}

The coefficients of the barycentric rational approximants
discussed in the previous subsections correspond to
singular vectors of the Loewner matrix $L\in\C^{n\times m}$.
This includes the case $n=m-1$ in Subsection~\ref{subsec:caseI};
a vector in the nullspace of $L$
can be understood as a singular vector corresponding to the
singular value~$\sigma_m=0$.

%We proceed to introduce auxiliary scaling matrices and
%the re-scaled Loewner matrix $\widehat{L}$.
The Loewner matrix is of the form~\eqref{eq.LoewnerMatrixform},
where the Cauchy matrix $C$ given in~\eqref{eq.defCauchyMat}
is a real matrix.
Thus, the Loewner matrix satisfies the representation~\eqref{eq.defA}
in Appendix~\ref{subsec.Asvd}, and we can define a re-scaled Loewner matrix,
\begin{subequations}
\begin{equation}\label{eq.defLhat}
 \widehat{L} = -\ii R L K\in\R^{n\times m},
\end{equation}
which is a real matrix as in~\eqref{eq.notatildeX}, and
\begin{equation}\label{eq.defRK0}
K=\diag{K_{11},\ldots,K_{mm}}\in\C^{m\times m},~~~\text{and}~~
R=\diag{R_{11},\ldots,R_{nn}}\in\C^{n\times n},
\end{equation}
are diagonal matrices as in~\eqref{eq.defLK}, with entries
\begin{equation}
\begin{aligned}
 K_{jj} =
\left\{ \begin{array}{ll}
(1-\ee^{-\ii y_j})/|1-\ee^{-\ii y_j}|,&~~~\ee^{-\ii y_j} \neq 1 ~~\text{and},\\
\ii,&~~~ \text{otherwise},
\end{array}\right.
\end{aligned}
\end{equation}
and
\begin{equation}\label{eq.defRs0}
\begin{aligned}
R_{kk} =
\left\{ \begin{array}{ll}
(1-\ee^{-\ii x_k})/|1-\ee^{-\ii x_k}|,&~~~\ee^{-\ii x_k} \neq 1 ~~\text{and},\\
\ii,&~~~ \text{otherwise}.
\end{array}\right.
\end{aligned}
\end{equation}
\end{subequations}
Following~\eqref{eq.notatildeXeq}, we also note the identity
\begin{equation}\label{eq.defLhatim}
\widehat{L} = 2 \imag(R C K^\ast),
\end{equation}
where $\imag(X)$ denotes the matrix of entry-wise imaginary parts of $X$.

The diagonal matrices $R$ and $K$ are unitary, and thus,
the matrices $L$ and $\widehat{L}$ are similar up to a complex phase and their singular values coincide.
Furthermore,
let the singular value decomposition of the matrix $\widehat{L}\in\R^{n\times m}$ be given by
\begin{equation}\label{eq.svdx1L}
\widehat{L} \widehat{V} = \widehat{U} S,
\end{equation}
where $\widehat{V}\in\R^{m\times m}$ and $ \widehat{U} \in\R^{n\times m}$
refer to orthonormal bases of right and left singular vectors,
respectively,
and $S=\diag{\sigma_1,\ldots,\sigma_m}\in\R^{m\times m}$
as in~\eqref{eq.svdL}.
We define
\begin{equation}\label{eq.wtildesvd}
\widetilde{w} = \ii K\widehat{V} e_m.
%~~~
%\text{with}~~\widehat{V}\in\R^{m\times m}~
%\text{as in~\eqref{eq.svdx1L}}.
\end{equation}
As a corollary of Proposition~\ref{prop.svdAmainprop}
in Appendix~\ref{subsec.Asvd}, the vector $\widetilde{w}$ is a singular
vector of $L$ corresponding to the singular value $\sigma_m$,
and as a consequence,
$$
\| L \widetilde{w} \|_2 = \sigma_m.
$$
Let $\widetilde{r}=\widetilde{n}/\widetilde{d}$ denote the
barycentric rational approximant with coefficients $\widetilde{w}$.
Then $\widetilde{r}$ minimizes the linearized error~\eqref{eq.linerrLWnorm}, since $\widetilde{w}$ is the singular vector corresponding to the smallest singular value, $\sigma_m$.
This also includes the case $n=m-1$ with $\sigma_m=0$, which is discussed in
Subsection~\ref{subsec:caseI}.

Let $\widetilde{w}=(\widetilde{w}_1,\ldots,\widetilde{w}_m)$
be given as in~\eqref{eq.wtildesvd}
and $f_j=\ee^{\ii y_j}$.
Then Proposition~\ref{prop.SVDidentities}
in Appendix~\ref{subsec.Asvd} implies
\begin{equation}\label{eq:fjwjiswjast0}
f_j \widetilde{w}_j = \widetilde{w}_j^\ast,~~~j=1,\ldots,m.
\end{equation}
As a consequence, the partial fractions $\widetilde{n}$ and $\widetilde{d}$ satisfy
\begin{equation}\label{eq:nisdast}
\widetilde{n}(x) = \sum_{j=1}^m \frac{f_j \widetilde{w}_j}{x-y_j}
 = \sum_{j=1}^m \frac{\widetilde{w}^\ast_j}{x-y_j} = \widetilde{d}(x)^\ast ,~~~\text{for $x\in\R$},
\end{equation}
and $\widetilde{r}$ has the representation
\begin{equation}\label{eq:rcayley}
\widetilde{r}(x)=\xi(x)^\ast \xi(x)^{-1},~~~\text{with}~~\xi(x)=\widetilde{d}(x).
\end{equation}
Let $\widetilde{r}$ be irreducible, then
this further implies that $\widetilde{r}$
has no poles on the real axis.
Furthermore, $\widetilde{r}$ is unitary, i.e., it satisfies
\begin{equation}
\tag{\ref{eq:unitarity}} 
\label{eq.rapproxtildeunitary}
|\widetilde{r}(x)|=1,~~~x\in\R.
\end{equation}

Concerning uniqueness of $w$ in~\eqref{eq.minLw},
we remark that singular
vectors are unique up to a complex phase for the \textit{non-degenerate} case,
$\sigma_{m-1} > \sigma_m$.
%, and thus, $w$ is only
%unique up to a complex phase in this case.
If the smallest singular value is degenerate,
any normalized linear combination of singular vectors
corresponding to $\sigma_m$ minimizes~\eqref{eq.minLw},
where in addition, each of these singular vector
has an arbitrary complex phase.

In the following proposition we consider barycentric rational approximants that minimize the linearized error but do not necessarily have coefficients~$\widetilde{w}$.

\begin{proposition}\label{prop.rminisuunitary}
Let $x_1,\ldots,x_{n}\in\R$ and $y_1,\ldots,y_m\in\R$ be given test and support nodes, respectively, with $n\geq m-1$.
Assume the smallest singular value of the corresponding Loewner matrix is non-degenerate,~i.e.,~$\sigma_{m-1}>\sigma_m$.
Let $r(x)=n(x)/d(x)$ denote the barycentric rational approximant with
coefficients $w_1,\ldots,w_m$
such that the linearized error~\eqref{eq.linerrLWnorm} is minimized.
Furthermore, assume~$r$ is irreducible.
Then $r$ has no poles on the real axis and satisfies
\begin{equation}
\tag{\ref{eq:unitarity}} 
|r(x)|=1,~~~x\in\R.
\end{equation}
\end{proposition}
\begin{proof}
The underlying vector of coefficients $w=(w_1,\ldots,w_m)^\top$
corresponds to a singular
vector of $L$ which is unique up to a complex phase
due to $\sigma_m$ being a non-degenerate.
Namely, we have $w=\ee^{\ii \phi}\widetilde{w}$ for $\widetilde{w}$
as given in~\eqref{eq.wtildesvd}
and some complex phase $\phi\in\R$.
Due to~\eqref{eq:fjwjiswjast0}, the entries of $w$ satisfy
$ f_j w_j = \ee^{2\ii\phi} w^\ast_j $.
Similar to~\eqref{eq:nisdast}, this implies $n(x)=\ee^{2\ii \phi} d(x)^\ast$ and we conclude $r(x)=\ee^{2\ii \phi} d(x)^\ast d(x)^{-1}$.
Thus, if $r$ is irreducible it has no poles on the real axis, and satisfies $|r(x)|=1$ for $x\in\R$.
\end{proof}

\begin{remark}\label{rmk:degandlincombination}
Proposition~\ref{prop.rminisuunitary} does not apply
if the smallest singular value of $L$ is degenerate,
i.e.,~$\sigma_m=\sigma_{m-1}$.
To clarify this remark we consider the following example.
Let $V$ denote the basis of right singular vectors of $L$,
then the vector $w=(\ee^{\ii \psi_1}V e_m+\ee^{\ii \psi_2}Ve_{m-1})/\sqrt{2}$
with arbitrary phases $ \psi_1, \psi_2\in\R $ minimizes $\|Lw\|_2$.
Thus, the barycentric rational approximant with coefficients $w$
minimizes the linearized error~\eqref{eq.linerrLWnorm}.
However, this approximant is not necessarily unitary.

On the other hand, the vector $\widetilde{w}$ given in~\eqref{eq.wtildesvd}
also satisfies~\eqref{eq:fjwjiswjast0} in case of $\sigma_m$ being
degenerate,
which entails unitarity of the respective barycentric rational approximant
even in the degenerate case.

Following Proposition~\ref{prop.SVDidentities}
in Appendix~\ref{subsec.Asvd}, $\widetilde{w}=\ii K\widehat{V}\zeta$
for any $\zeta\in\R^m$
satisfies $f_j\widetilde{w}_j = \widetilde{w}_j^\ast$.
Thus, coefficients corresponding to this vector also generate
a unitary rational approximant. However, such an approximant
does not minimize the linearized error.
Similar results hold for $\widetilde{w}=\ii K\widehat{V}\zeta$ with 
$\zeta\in\ii\R^m$.
\end{remark}

\subsection{Minimizing a weighted linearized error}
\label{subsec:caseIImu}

We proceed to generalize results of the previous subsections by
considering barycentric rational approximants of the form~\eqref{eq.defr}
that minimize a weighted linearized error over the test nodes.
In the unweighted case, the deviation in the approximation $r(x)\approx \ee^{\ii x}$ at a test node $x_k$ is $r(x_k)-\ee^{\ii x_k}$. 
In the weighted case, we are given weights,
\begin{align}\label{eq:defweightedgam}
\mu_1,\ldots,\mu_n>0,
\end{align} 
and consider a weighted deviation at test nodes,
\begin{equation}\label{eq.weightederr0}
\mu_k^{1/2} ( r(x_k) - \ee^{\ii x_k} ).
\end{equation}
We slightly modify the Loewner matrix~\eqref{eq.defLentries}
to include the weights $\mu_k$. Namely, we define
$$
L_{\mu}\in\C^{n\times m}~~~\text{with}~~
(L_{\mu})_{kj} =
\frac{\mu_k^{1/2}( \ee^{\ii x_k} - \ee^{\ii y_j})}{x_k-y_j}.
$$
Similar to~\eqref{eq.LoewnerMatrixform},
this matrix satisfies the representation
\begin{equation}\label{eq.weightedLoewner}
L_{\mu}
= S_F M - M S_f,~~~ \text{with}~~M=S_\mu C\in\R^{n\times m},
\end{equation}
where $C$ denotes the Cauchy matrix as given in~\eqref{eq.defCauchyMat}
and $S_\mu=\diag{\mu_1^{1/2},\ldots,\mu_n^{1/2}}\in\R^{n\times n}$.
Similar to~\eqref{eq.Lwtoxj},
applying a vector of coefficients $w=(w_1,\ldots,w_m)^\top\in\C^m$
to the weighted Loewner matrix yields
$$
(L_{\mu} w)_k = \mu_k^{1/2}( \ee^{\ii x_k} d(x_k) - n(x_k)),
$$
where $d$ and $n$ are partial fractions of the respective
barycentric rational representation as in~\eqref{eq.defr}.
Thus, we have
\begin{equation}\label{eq.weightedlinerr}
\|L_{\mu} w \|_2 = \left( \sum_{k=1}^n \mu_k | \ee^{\ii x_k} d(x_k) - n(x_k) |^2\right)^{1/2}.
\end{equation}
This corresponds to a weighted linearized error
which is related to the weighted deviation~\eqref{eq.weightederr0}.

Considering a barycentric rational approximant $r(x)=n(x)/d(x)$ which minimizes~\eqref{eq.weightedlinerr},
we let coefficients $w_1,\ldots,w_m\in \C$ of $r$
correspond to a right singular vector of~$L_\mu$, namely, the right singular vector
corresponding to the smallest singular value of~$L_\mu$.
The matrices $L$ and $L_\mu$ both satisfy the representation~\eqref{eq.defA} in Appendix~\ref{subsec.Asvd}.
Thus, the results of the previous subsections  -- in particular, Proposition~\ref{prop.rminisuunitary}
and the Cayley-type representation given in~\eqref{eq:rcayley} -- also apply
when considering a weighted linearized error.

\section{Non-interpolatory rational approximation}
\label{sec:caseIII}

In the present section we consider a barycentric rational approximant
$r_b$ based on the representation~\eqref{eq.defr2}
that satisfies some accuracy conditions at $n$ test nodes, where $n>m-1$.
We proceed with a setting similar to Section~\ref{sec:LoewnerMat}:
Let the degree $m$ be fixed, and let $y_1,\ldots,y_m\in\R$
and $x_1,\ldots,x_n$ be given support and test nodes, respectively.
In contrast to the barycentric rational approximant $r$
based on the representation~\eqref{eq.defr}, which is the topic of the previous
section, the approximant $r_b$
does not necessarily interpolate at the support nodes.
On the other hand, $r_b$ has
more coefficients which are free to be chosen such that
a higher accuracy at the test nodes can be achieved
for this representation.
We aim to choose the coefficients of $r_b$ such that
a weighted linearized error is minimized in the progress.
This approach is motivated by the AAA--Lawson algorithm~\cite{NT20}.
The AAA--Lawson algorithm runs in two phases.
First, the AAA method is applied, which results in an approximant
$r$ of the form~\eqref{eq.defr}, and support nodes $y_1,\ldots,y_m$.
In the second phase, these support nodes are
put in the representation~$r_b$ in~\eqref{eq.defr2},
and the coefficients $\alpha_1,\ldots,\alpha_m$ and $\beta_1,\ldots,\beta_m$
are determined such that a weighted linearized error is minimized.
This phase runs iteratively, adapting the weights
of the underlying linearized error with each iteration.
A similar approach is already suggested in~\cite{NST18}
and appears in~\cite[Section~8]{FNTB18} for the approximation of real functions.
Focusing on the unitarity of the generated approximant,
we assume the support nodes $y_1,\ldots,y_m$ and the weights to be given.

We recall the barycentric representation
given in~\eqref{eq.defr2},
\begin{equation}\label{eq:ralpasnd}
r_{b}(x) = \sum_{j=1}^m \frac{\alpha_j}{x-y_j}
\Big/ \sum_{j=1}^m \frac{\beta_j}{x-y_j} =: n_{b}(x) / d_{b}(x),
\end{equation}
where $\alpha_1,\ldots,\alpha_m\in\C$ and
$\beta_1,\ldots,\beta_m\in\C$ denote coefficients which have to be chosen
to minimize the linearized error.
Let
$$
\mu_1,\ldots,\mu_n>0
$$
denote given weights
which represent scaling factors of the deviation in the approximation $r(x) \approx \ee^{\ii x}$ at the test nodes,
$$
\mu_k^{1/2} ( r_{b}(x_k) - \ee^{\ii x_k} ).
$$
In the current setting, $r_{b}$ does not necessarily interpolate $\ee^{\ii x}$ at the support nodes,
and it is reasonable to include support nodes in the set of test nodes
such that accuracy at the support nodes can also be enforced.
Thus, in the present section
we allow test nodes to coincide with support nodes.

To derive a weighted linearized error
for the barycentric rational approximant
$r_b=n_b/d_b$ in~\eqref{eq:ralpasnd},
we first simplify $r_{b}(x) \approx \ee^{\ii x}$
to $n_{b}(x) \approx \ee^{\ii x} d_{b}(x) $.
Evaluated at a test node $x_k$ which is not a support node,
we consider the linearized deviation
\begin{subequations}\label{eq:linerrorlawson0}
\begin{equation}\label{eq:linerrorlawson1}
\mu_k^{1/2} ( n_b(x_k) - \ee^{\ii x_k} d_b(x_k)).
%\Big| \ee^{\ii x_k} \sum_{j=1}^m \frac{\beta_j}{x_k-y_j}
%- \sum_{j=1}^m \frac{\alpha_j}{x_k-y_j} \Big|.
\end{equation}
If the test node $x_k$ is also a support node,
there is an index~$j_k$ such that~$x_k=y_{j_k}$, and the representation $r_b$ as given in~\eqref{eq:ralpasnd}
can not be evaluated at $x=x_k$ in a direct manner due to the presence of the partial fractions $\frac{\alpha_{j_k}}{x-y_{j_k}}$ and $\frac{\beta_{j_k}}{x-y_{j_k}}$.
However, in the limit $x\to x_k$ we attain
$$
\lim_{x \to x_k} r_b(x) = \alpha_{j_k}/\beta_{j_k}, %~~~\text{for}~~x_k=y_{j_k},
$$
and~\eqref{eq:linerrorlawson1} is replaced by
\begin{equation}\label{eq:linerrorlawson2}
\mu_k^{1/2} ( \alpha_{j_k} - \ee^{\ii x_k} \beta_{j_k} ),~~~\text{for}~~x_k=y_{j_k}.
\end{equation}
\end{subequations}
The deviations in~\eqref{eq:linerrorlawson0}
can be represented as a matrix-vector product
of an expanded Loewner matrix with the coefficient vectors
$\alpha=(\alpha_1,\ldots,\alpha_m)^\top\in\C^m$ and
$\beta=(\beta_1,\ldots,\beta_m)^\top\in\C^m$.
We introduce a modified Cauchy matrix $C'\in\R^{n\times m}$ as follows.
\begin{itemize}
\item Let $k$ be an index with $x_k\notin\{y_1,\ldots,y_m\}$.
Then we define the $k$-th row of $C'$
analogously to the $k$-th row of the Cauchy matrix $C$ in~\eqref{eq.defCauchyandDiag}.
This yields
\begin{subequations}\label{eq:CalpCbet}
\begin{equation}\label{eq:CalpCbet1}
(C' \alpha)_k = \sum_{j=1}^m \frac{\alpha_j}{x_k-y_j} =  n_{b}(x_k),~~~\text{and}~~
(C' \beta)_k =  \sum_{j=1}^m \frac{\beta_j}{x_k-y_j} = d_{b}(x_k).
\end{equation}
\item Otherwise, for an index $k$ with $x_k\in\{y_1,\ldots,y_m\}$
we have an index $j_k$ such that $x_k=y_{j_k}$,
and we define the $k$-th row of $C'$ by
$$
C'_{kj} = \left\{
\begin{array}{ll}
1,&~~~\text{for $j=j_k$, and}\\
0,&~~~\text{for $j\in\{1,\ldots,m\}\setminus\{j_k\}$},
\end{array}
\right.~~~\text{for $k,j_k$ with $x_k=y_{j_k}$.}
$$
For this index $k$, the matrix-vector products in~\eqref{eq:CalpCbet1} yield
\begin{equation}\label{eq:CalpCbet2}
(C' \alpha)_k = \alpha_{j_k},~~~\text{and}~~
(C' \beta)_k = \beta_{j_k}.
\end{equation}
\end{subequations}
\end{itemize}

Furthermore, let $S_F\in\C^{n\times n}$ denote the diagonal matrix with diagonal entries $\ee^{\ii x_k}$
as in~\eqref{eq.defSfSfap}, and
let $M=S_\mu C'\in\R^{n\times m}$
where $S_\mu=\diag{\mu_1^{1/2},\ldots,\mu_n^{1/2}}\in\R^{n\times n}$.
The matrix-vector products $M\alpha$ and $M\beta$
can be evaluated similar to~\eqref{eq:CalpCbet}.
For the concatenated vector $\gamma=[\alpha;\beta]\in\C^{2m}$
we have
\begin{subequations}\label{eq.linerrorlawson0}
\begin{equation}\label{eq.linerrorlawsonmatform}
([M | -S_F M]\gamma)_k =
\left\{\begin{array}{ll}
\mu_k^{1/2}( n_b(x_k) - \ee^{\ii x_k} d_b(x_k),&~~
x_k\notin\{y_1,\ldots,y_m\}~~\text{and},\\
\mu_{k}^{1/2} (\alpha_{j_k}- \ee^{\ii x_k} \beta_{j_k}),&~~
x_k=y_{j_k}~~\text{for some $j_k$},
\end{array}\right.
\end{equation}
where $[M| -S_F M]$ is to be understood as the $n\times 2m$ complex matrix
obtained by concatenating the $n \times m$ matrices $M$ and $-S_F M$.
The Euclidean norm of $ [M | -S_F M]\gamma $ quantifies the deviations~\eqref{eq:linerrorlawson0}
over the test nodes (which may include support nodes) and is also referred to as linearized error
in the present section.

In the following,
we consider barycentric rational approximants $r_b$
that attain a minimal linearized error.
Namely, the vector $\gamma=[\alpha;\beta]\in\C^{2m}$,
where $\alpha$ and $\beta$ refer to the coefficients of~$r_b$, satisfies
\begin{equation}\label{eq.minlawson}
\gamma=\argmin_{v\in\C^{2m},~\|v\|_2=1} \big\|[M| -S_F M]v \big\|_2.
\end{equation}
\end{subequations}
Similar to the coefficient vector $w$ in~\eqref{eq.minLw},
the vector $\gamma$ in~\eqref{eq.minlawson} is accessible using a singular value
decomposition of the matrix $[M| -S_F M]\in\C^{n\times 2m}$.

We proceed with some auxiliary properties of right singular vectors
of $[M| -S_F M]$.
These results are given in detail in Appendix~\ref{subsec.appendixcaseIII}.
Let $R=\diag{R_{11},\ldots,R_{nn}}\in\C^{n\times n}$ be given as in~\eqref{eq.defLK}
in Appendix~\ref{subsec.Asvd},~i.e.,
\begin{equation}\label{eq:defRkkforBhat}
\begin{aligned}
R_{kk} =
\left\{ \begin{array}{ll}
(1-\ee^{-\ii x_k})/|1-\ee^{-\ii x_k}|,&~~~\ee^{-\ii x_k} \neq 1 ~~\text{and},\\
~ \ii,&~~~ \text{otherwise},
\end{array}\right.
\end{aligned}
\end{equation}
and let
\begin{equation}\label{eq.defBhat}
\widehat{B} = [\real(R) M| -\imag(R) M]\in\R^{n\times 2m},
\end{equation}
where $\real(R)$ and $\imag(R)$ refer to the matrices
of entry-wise real and imaginary parts of $R$, respectively.
The matrix $\widehat{B}$ satisfies~\eqref{eq:defBhat}
in Appendix~\ref{subsec.appendixcaseIII} with
$\Theta=RM$.
Let $\widehat{V}_B\in\R^{2m\times 2m}$ denote a real
orthonormal basis of right singular vectors of $\widehat{B}$
as in~\eqref{eq:Bhatsvd},
and let
\begin{subequations}\label{eq.defgam}
\begin{equation}\label{eq.defhatgam}
\widehat{\gamma}=\widehat{V}_B e_{2m}\in\R^{2m},~~~
\widehat{\gamma}=[\widehat{\alpha};\widehat{\beta}],
\end{equation}
denote a right singular vector of $\widehat{B}$ corresponding
to the smallest singular value of $\widehat{B}$.
We introduce the vector
\begin{equation}\label{eq.deftildegam}
\widetilde{\gamma} =
\frac{1}{\sqrt{2}}
\left(\begin{array}{c}
\widehat{\alpha} + \ii\widehat{\beta}\\
\widehat{\alpha} - \ii\widehat{\beta}
\end{array}\right)\in\C^{2m}, ~~~~
\widetilde{\gamma}=[\widetilde{\alpha};\widetilde{\beta}],
\end{equation}
which corresponds to the vector $\widetilde{\gamma}$ in~\eqref{eq.deftildegamx}
in Appendix~\ref{subsec.appendixcaseIII}.
Following Corollary~\ref{cor:svabsabisequal} in Appendix~\ref{subsec.appendixcaseIII},
$\widetilde{\gamma}$ satisfies~\eqref{eq.minlawson},
and we have the identity
\begin{equation}\label{eq.alpisbetast1}
\widetilde{\alpha}_j = \widetilde{\beta}^\ast_j,~~~j=1,\ldots,m.
\end{equation}
\end{subequations}
Thus, the barycentric rational representation
$\widetilde{r}_b=\widetilde{n}_b/\widetilde{d}_b$
with coefficients $\widetilde{\alpha}$ and $\widetilde{\beta}$
minimizes the linearized error, and $\widetilde{n}_b$ and $\widetilde{d}_b$ satisfy
\begin{equation}\label{eq.nisdastalph1}
\widetilde{n}_b(x)
=\sum_{j=1}^m \frac{\widetilde{\alpha}_j}{x-y_j}
 = \sum_{j=1}^m \frac{\widetilde{\beta}^\ast_j}{x-y_j}
=\widetilde{d}_b(x)^\ast,~~~\text{for $x\in\R$}.
\end{equation}
As a consequence, we have the Cayley-type representation
\begin{equation}\label{eq:rbcayley}
\widetilde{r}_b(x) = \xi(x)^\ast \xi(x)^{-1},~~~\text{with}~~\xi(x) = \widetilde{d}_b(x),
\end{equation}
and for the case that $\widetilde{r}_b$ is irreducible,
this implies that $\widetilde{r}_b$ has no poles on the real axis and
\begin{equation}
\tag{\ref{eq:unitarity}}
|\widetilde{r}_b(x)|=1,~~~x\in\R.
\end{equation}

If the smallest singular value of $[M| -S_F M]$ is non-degenerate,
then a similar result carries over to any barycentric rational approximant that minimizes the respective linearized error.

\begin{proposition}\label{prop:maincaseIII}
%Let support nodes $y_1,\ldots,y_m\in\R$
%and weigts $\mu_1,\ldots,\mu_m>0$ be given.
Let the smallest singular value of $[M| -S_F M]$ be non-degenerate,
and let $\gamma=[\alpha;\beta]$ satisfy~\eqref{eq.minlawson}.
Assume the generated barycentric rational approximant $r_b$ is irreducible.
Then $r_b$ has no poles on the real axis and satisfies
\begin{equation}
\tag{\ref{eq:unitarity}}
|r_b(x)|=1,~~~x\in\R.
\end{equation}
\end{proposition}
\begin{proof}
For the case that the smallest singular value of $[M| -S_F M]$ is non-degenerate, coefficients $\alpha$ and $\beta$ that satisfy~\eqref{eq.minlawson}
correspond to $\alpha=\ee^{\ii \phi} \widetilde{\alpha}$ and
$\beta=\ee^{\ii \phi} \widetilde{\beta}$ for a phase $\phi\in\R$,
and $\widetilde{\alpha}$ and $\widetilde{\beta}$ as given in~\eqref{eq.deftildegam}.
Due to~\eqref{eq.alpisbetast1}, we have $\alpha_j=\ee^{2\ii \phi} \beta^\ast_j$ for $j=1,\ldots,m$.
Similar to~\eqref{eq.nisdastalph1}, this yields
 $n_b(x)=\ee^{2\ii \phi} d_b(x)^\ast$ for $n_b$ and $d_b$ given in~\eqref{eq:ralpasnd}.
Thus, the zeros of $d_b$ and $n_b$ located on the real axis coincide.
Together with the assumption that $r_b=n_b/d_b$ is irreducible, this entails that $r_b$ has no poles
on the real axis.
Furthermore, we have $|r_b(x)|=1$ for $x\in\R$.
\end{proof}

Considering the restriction to a non-degenerate case in Proposition~\ref{prop:maincaseIII},
we also refer to Remark~\ref{rmk:degandlincombination} in the previous section.

\begin{remark}\label{rmk:rbdeg}
Analogously to the barycentric rational approximant $\widetilde{r}_b$
with coefficients corresponding to $\widetilde{\gamma}$ as given in~\eqref{eq.deftildegam},
a barycentric rational approximant with coefficients
corresponding to $\widetilde{\gamma}' = Q\widehat{V}_B \zeta\in\C^{2m}$
(with $Q$ and $\widehat{V}_B$ as in Proposition~\ref{prop:svdMSFMgen}
in Appendix~\ref{subsec.appendixcaseIII})
with $\zeta\in\R^{2m}$ or $\zeta\in\ii\R^{2m}$ is unitary as well.
This also holds true when singular values of $[M|-S_F M]$ are degenerate.
\end{remark}

\section{Advantages of Cayley-type representations}
\label{subsection:Cayleyform}

In the present section we give some remarks concerning the computation
of the interpolatory and non-interpolatory approximants introduced in the previous sections.
We have shown the unitarity property~\eqref{eq:unitarity} for approximants which minimize
linearized errors in the previous sections.
The coefficients of the respective barycentric rational approximants
are unique up to a complex phase.
In practice, we suggest computing coefficients using the re-scaled Loewner matrix~$\widehat{L}\in\R^{n\times m}$~\eqref{eq.defLhat} instead of the original Loewner matrix~$L\in\C^{n\times m}$ in the context of the AAA method, or the modified Loewner-type matrix~$\widehat{B}\in\R^{n\times 2m}$~\eqref{eq.defBhat} instead of~$[M|-S_F M]\in\C^{n\times 2m}$ in the context of the AAA--Lawson method.
Furthermore, we suggest utilizing a Cayley-type representation of the rational approximant.
These modifications can be implemented for the AAA and AAA--Lawson methods~\cite{NST18,NT20}, as illustrated in Algorithm~\ref{alg:modAAA} and~\ref{alg:modAAALawson} below, which gives some advantages in terms of computational cost and numerical stability of the unitarity property.
%We suggest computing coefficients using the re-scaled Loewner matrix $\widehat{L}\in\R^{n\times m}$~\eqref{eq.defLhat} or the modified Loewner-type matrix $\widehat{B}\in\R^{n\times 2m}$~\eqref{eq.defBhat} in practice.
%This yields some advantages compared to the original AAA and AAA--Lawson algorithms in~\cite{NST18,NT20} which make use of the Loewner matrix $L\in\C^{n\times m}$ or an expanded Loewner matrix $[M|-S_F M]\in\C^{n\times 2m}$, respectively.

We first consider the interpolatory case
and the approximant $\widetilde{r}$~\eqref{eq:rcayley}
of the representation~\eqref{eq.defr}
with coefficients $\widetilde{w}_1,\ldots,\widetilde{w}_m$
as given in~\eqref{eq.wtildesvd}.
\begin{enumerate}
\item[(i)] The coefficients $\widetilde{w}_1,\ldots,\widetilde{w}_m$
can be computed via a singular value decomposition
of the re-scaled Loewner matrix
$\widehat{L}$ using real arithmetic.
This typically reduces computational cost
compared with the original AAA algorithm, which uses
a singular value decomposition of a complex Loewner matrix.
\item[(ii)] Exploiting the Cayley-type representation $\widetilde{r}=\xi^\ast \xi^{-1}$
reduces computational cost when $\widetilde{r}$ has to be evaluated, or
zeros and poles of $\widetilde{r}$ are required.
In particular, the zeros of $\widetilde{r}$ are the complex conjugate of its poles.
\item[(iii)]
The complex phases of the coefficients $\widetilde{w}_1,\ldots,\widetilde{w}_m$
originate from applying the complex diagonal matrix $\ii K$ in~\eqref{eq.wtildesvd}
on the real vector $\widehat{V} e_m$.
Thus, these complex phases are exact up to machine precision, and are not affected by errors that occur from the underlying singular value decomposition.
As a consequence, the identity $f_j w_j= w_j^\ast$~\eqref{eq:fjwjiswjast0}
is true up to machine precision.
This results in an improved numerical stability on the
unitarity of $\widetilde{r}$ compared to~$r$. For the deviation of the unitarity property
of approximants generated by the original and modified AAA algorithms in computer arithmetic, see~\figref{fig.UnitarityRoundoff}.
\item[(iv)] In a degenerate or close to degenerate case, the approximant $\widetilde{r}$
with coefficients $\widetilde{w}_1,\ldots,\widetilde{w}_m$
remains unitary, see also Remark~\ref{rmk:degandlincombination}.
\end{enumerate}
In Algorithm~\ref{alg:modAAA}, we illustrate the modifications (i) and (ii) for the AAA method, which yields a modified AAA algorithm utilizing (i)--(iv).
This algorithm is based on a simplified version of \texttt{aaa.m}
in the \texttt{chebfun} package~\cite{DHT14}.
%Major modifications are highlighted by `$\leftarrow$' symbols in the code.
Note that only a minimal version of the original AAA algorithm -- adequate for illustrating the modification based on the Cayley-type representation --
is described in the present work. For further details of AAA we refer to~\cite{NST18},
and the short summary given in Section~\ref{sec:AAA}.
%The auxiliary function \texttt{reval\_expi}
%and \texttt{prz\_expi} refer to functions similar to \texttt{reval.m}
%and \texttt{prz.m} in the \texttt{chebfun} package, respectively.
%However, \texttt{reval\_expi} and \texttt{prz\_expi} utilize~(ii) above in a direct manner.

\begin{algorithm}
\caption{Pseudocode for a modified AAA algorithm
generating a barycentric rational approximant of $\ee^{\ii x}$.
See~\texttt{aaa.m} in the~\texttt{chebfun} package and
\cite[Fig.~4.1]{NST18} for the original code in full detail.}
\label{alg:modAAA}
\begin{algorithmic}[1]
\Require{$(x_1,\ldots,x_N),m_{\text{max}},n_{\text{lawson}}$}
\State{$F \Leftarrow (\ee^{\ii x_1},\ldots,\ee^{\ii x_N})$}
\State{$R \Leftarrow \text{diag}\big({(1-\ee^{-\ii x_1})}/{|1-\ee^{-\ii x_1}|},\ldots,{(1-\ee^{-\ii x_N})}/{|1-\ee^{-\ii x_N}|} \big)$}
\State{for all $j$ with $(\ee^{-\ii x_j}==1)$ set $R_{jj}= \ii$}
\Comment{the matrix $R$ as in~\eqref{eq.defRs0}}
\State{$r \Leftarrow \text{mean}(F) \cdot(1,\ldots,1)\in\mathbb{C}^N$}
\For{$m = 1,\ldots,m_{\text{max}}$}
\State{$ j \Leftarrow \argmax ( |F - r| )$}
\State{$y_m \Leftarrow x_j$}
\State{$ f_m \Leftarrow F_j$}
\State{$ k_m \Leftarrow R_{jj}$}
\State{remove $x_j$ from the list of test nodes, i.e.,~$x\Leftarrow (x_1,\ldots,x_{j-1},x_{j+1},\ldots,x_{N-m+1})$

and remove the respective columns of $C$ and $F$ and entries of $R$}
\State{$c_{\text{new}} \Leftarrow \big(\frac{1}{x_1-y_m},\ldots,\frac{1}{x_{N-m}-y_m} \big)^\top$}
\State{$    C \Leftarrow [C ~|~ c_{\text{new}}]$}
\State{$    K \Leftarrow \text{diag}(k_1,\ldots,k_m)$}
\State{$  \widehat{L} \Leftarrow 2\imag(R \cdot C \cdot \conj(K))\in\mathbb{R}^{(N-m)\times m}$}
\Comment{the re-scaled Loewner matrix~\eqref{eq.defLhatim}.}
\State{$    (\widehat{U}, S, \widehat{V}) \Leftarrow \texttt{svd}(\widehat{L})$}
\State{$    \widetilde{w} \Leftarrow \ii K \widehat{V} e_m \in \mathbb{C}^{m}$}
\Comment{see~\eqref{eq.wtildesvd}}
\State{$  \xi  \Leftarrow C\cdot \widetilde{w} $}
\State{$r  \Leftarrow ( \xi^\ast_1/\xi_1,\ldots,  \xi^\ast_{N-m}/\xi_{N-m})$}
\Comment{$=(r(x_1),\ldots,r(x_{N-m}))$
using $r=\xi^\ast/\xi$}
\EndFor
\If{$n_{\text{lawson}}>0$}
\Comment{run Algorithm~\ref{alg:modAAALawson}}
\State{$   \widetilde{w} \Leftarrow  \text{minimax\_lawson}(x=(x_1,\ldots,x_{N-m}),y=(y_1,\ldots,y_m),n_{\text{lawson}})$}
\EndIf
\State{\Return{$ x \mapsto r(x) = \sum_{j=1}^m \frac{\widetilde{w}^\ast_j}{x-y_j} / \sum_{j=1}^m \frac{\widetilde{w}_j}{x-y_j}$}}
\Comment{using the Cayley-type representation~\eqref{eq:rcayley}

\Comment{for the case $n_{\text{lawson}}>0$ with
$\widetilde{w}_j=\widetilde{\beta}_j$, the generated $r(x)$ corresponds to~\eqref{eq:rbcayley}}}
\end{algorithmic}
\end{algorithm}

We proceed to discuss advantages of
utilizing the non-interpolatory approximant $\widetilde{r}_b$
in~\eqref{eq:rbcayley} with the coefficients
$\widetilde{\alpha}_1,\ldots,\widetilde{\alpha}_m$
and $\widetilde{\beta}_1,\ldots,\widetilde{\beta}_m$
given in~\eqref{eq.defgam}.
The following points are similar to~(i)-(iv) above.
\begin{enumerate}
\item[(i-b)] The coefficients $\widetilde{\alpha}_1,\ldots,\widetilde{\alpha}_m$
and $\widetilde{\beta}_1,\ldots,\widetilde{\beta}_m$
can be computed via a singular value decomposition
of $\widehat{B}\in\R^{n\times 2m}$~\eqref{eq.defBhat}
using real arithmetic.
This reduces computational cost compared with the original AAA--Lawson algorithm which uses the singular value decomposition of the complex matrix $[M|-S_F M]\in\C^{n\times 2m}$.
\item[(ii-b)] Exploiting the Cayley-type representation $\widetilde{r}_b=\xi^\ast \xi^{-1}$
given in~\eqref{eq:rbcayley}
reduces computational cost when $\widetilde{r}_b$ has to be evaluated, or
zeros and poles of $\widetilde{r}_b$ are required.
\item[(iii-b)] Using the formula~\eqref{eq.defgam}
to compute the coefficients $\widetilde{\alpha}_1,\ldots,\widetilde{\alpha}_m$
and $\widetilde{\beta}_1,\ldots,\widetilde{\beta}_m$
helps avoid errors on the complex phases of the coefficients,
which can occur during the computation of a singular value decomposition of $[M|-S_F M]$.
This also results in an improved numerical stability on the
unitarity of $\widetilde{r}_b$ compared to $r_b$. For the deviation of the unitarity property
of approximants generated by the original and modified AAA--Lawson algorithms in computer arithmetic, see~\figref{fig.UnitarityRoundoff}.
\item[(iv-b)] In a degenerate or close to degenerate case, the approximant $\widetilde{r}_b$
with coefficients $\widetilde{\alpha}_1,\ldots,\widetilde{\alpha}_m$
and $\widetilde{\beta}_1,\ldots,\widetilde{\beta}_m$
remains unitary, see also Remark~\ref{rmk:rbdeg}.
\end{enumerate}
In Algorithm~\ref{alg:modAAALawson}, we illustrate the modifications (i-b) and (ii-b) for the \textit{minimax} iteration which is part of the AAA--Lawson method. This algorithm is based on a simplified version of \texttt{aaa.m} in the \texttt{chebfun} package~\cite{DHT14}. In combination with Algorithm~\ref{alg:modAAA} (and `$n_{\text{lawson}}>0$'), Algorithm~\ref{alg:modAAALawson} yields a modified version of the AAA--Lawson algorithm utilizing (i-b)--(iv-b).
% and
%for the non-interpolatory approximation that minimizes the linearized error.
%Main modifications in Algorithm~\ref{alg:modAAA} and~\ref{alg:modAAALawson}
%introduced in the present work are highlighted by `$\leftarrow$' symbols.

\begin{algorithm}
\caption{Pseudocode for a
modified \textit{minimax} iteration which is called by Algorithm~\ref{alg:modAAA}
to generate coefficients of a near-best rational approximant
of $\ee^{\ii x}$.
This algorithm is based on \texttt{aaa.m} in the \texttt{chebfun} package
and yields a modified version of the AAA--Lawson method.}
\label{alg:modAAALawson}
\begin{algorithmic}[1]
\Require{$x=(x_1,\ldots,x_n),y=(y_1,\ldots,y_m),n_{\text{lawson}}$}
\State{$N\Leftarrow n+m$}
\State{$x\Leftarrow (x_1,\ldots,x_n,y_1,\ldots,y_m)\in\mathbb{R}^N$}
\Comment{include support nodes to the set of test nodes}
\State{$\mu \Leftarrow (1,\ldots,1)\in\mathbb{R}^N$}
\State{$F \Leftarrow (\ee^{\ii x_1},\ldots,\ee^{\ii x_N})$}
\State{$R \Leftarrow \text{diag}\big({(1-\ee^{-\ii x_1})}/{|1-\ee^{-\ii x_1}|},\ldots,{(1-\ee^{-\ii x_N})}/{|1-\ee^{-\ii x_N}|} \big)$}
\State{for all $j$ with $(\ee^{-\ii x_j}==1)$ set $R_{jj}= \ii$}
\Comment{the matrix $R$ as in~\eqref{eq:defRkkforBhat}}
\State{$C\Leftarrow0\in\mathbb{C}^{N\times m}$}
\Comment{the Cauchy matrix $C'$, some test and support nodes coincide}
\For{$j = 1,\ldots,n$}
\State{$C_{j, :} \Leftarrow \big(\frac{1}{x_j-y_1},\ldots,\frac{1}{x_{j}-y_m} \big) $}
\EndFor
\For{$j = 1,\ldots,m$}
\State{$C_{n+j,j} \Leftarrow 1$}
\EndFor
\State{$A \Leftarrow [\real(R)\cdot C| -\imag(R) \cdot C]\in\mathbb{R}^{N\times 2m}$}
\For{$steps = 1,\ldots,n_{\text{lawson}}$}
\State{$S_\mu \Leftarrow \text{diag}( \sqrt{\mu_1},\ldots,\sqrt{\mu_N})$}
\State{$(\widehat{U},\widehat{S},\widehat{V}) \Leftarrow \texttt{svd}(S_\mu\cdot A)$}
\Comment{\texttt{svd} of $\widehat{B}=S_\mu\cdot A$, see~\eqref{eq.defBhat}}
\State{$\widehat{\gamma} \Leftarrow \widehat{V} e_{2m} \in\mathbb{R}^{2m}$}
\Comment{see~\eqref{eq.defhatgam}}
\State{$\widetilde{\beta} = (\widetilde{\beta}_1,\ldots,\widetilde{\beta}_m)$ with $\widetilde{\beta}_j \Leftarrow (\widehat{\gamma}_j - \ii\widehat{\gamma}_{j+m})/\sqrt{2}$}
\Comment{see~\eqref{eq.deftildegam}}
\State{$\xi \Leftarrow 0\in\mathbb{C}^{N}$}
\For{$\ell = 1,\ldots,n$}
\State{$\xi_\ell \Leftarrow \sum_{j=1}^m \frac{\widetilde{\beta}_j}{x_{\ell} - y_j}$}
\EndFor
\For{$\ell = 1,\ldots,m$}
\State{$\xi_{n+\ell} \Leftarrow \widetilde{\beta}_\ell$}
\EndFor
\State{$r \Leftarrow (\xi^\ast_1/\xi_1,\ldots,\xi^\ast_N/\xi_N)$}
\Comment{using the Cayley-type representation~\eqref{eq:rbcayley}}
\State{$\varepsilon \Leftarrow F - r$}
\State{$\mu \Leftarrow (\mu_1 |\varepsilon_1|,\ldots,\mu_N |\varepsilon_N|)$}
\State{$\mu \Leftarrow \mu/\|\mu\|_{\infty}$}
\EndFor
\State{\Return $\widetilde{\beta}=(\widetilde{\beta}_1,\ldots,\widetilde{\beta}_m)$}
\end{algorithmic}
\end{algorithm}

%%%%%%%%%%%%%%%%%%%%%
\section{Unitarity of AAA and AAA--Lawson methods}\label{sec:AAA}

We briefly sketch the AAA and AAA--Lawson algorithms
to illustrate that Proposition~\ref{prop.rminisuunitary} and~\ref{prop:maincaseIII}
apply in these settings, respectively.
This shows unitarity~\eqref{eq:unitarity} of the generated approximants.

\subsection{Unitarity of the AAA method}
\label{subsec:AAA}

Let us consider the application of the AAA algorithm~\cite{NST18} to the approximation
of $\ee^{\ii x}$. This algorithm aims to generate a
barycentric rational approximant $r$ using the representation~\eqref{eq.defr}
which is accurate at a given set of test nodes $x_1,\ldots,x_n\in\R$,
while the degree $m$ of the approximant is increased in an iterative manner.
%This method can also provide approximants which are accurate on a given interval
%by disretizing the interval in an adaptive manner.
The initial iteration step uses an approximant
with a single support node, i.e., $m=1$,
and, the iteration proceeds as follows.
\begin{enumerate}
\item[(a)] The number of support nodes,~i.e., the degree $m$, is increased.
In particular, the test node for which the previously computed approximant
yields the largest deviation
from $\ee^{\ii x}$ is added to the set of support nodes, and removed as a test node.
\item[(b)] The coefficients $w_1,\ldots,w_m\in\C$
are computed such that the linearized error
is minimized~\eqref{eq.minLw}, which generates a new approximant.
\item[(c)] Accuracy conditions are tested for the new approximant.
If error tolerances are satisfied, AAA returns the generated approximant.
Otherwise, the iteration proceeds with (a).
\end{enumerate}
In addition, \textit{Froissart doublets} can be detected by AAA.
If Froissart doublets occur,
specific nodes are removed from the set of
support and test nodes, and the coefficients are re-computed by minimizing
the linearized error for the new sets of nodes.

Due to the fact that coefficients $w_1,\ldots,w_m$
minimize the linearized error as in~\eqref{eq.minLw},
Proposition~\ref{prop.rminisuunitary} applies
and shows that the generated approximant has no poles on the real axis
and is unitary.
This requires further conditions of Proposition~\ref{prop.rminisuunitary}
to hold true. Namely, the smallest singular value of
the underlying Loewner matrix has to be non-degenerate,
and the generated approximant has to be irreducible.

\subsection{Unitarity of the AAA--Lawson method}
\label{subsec:AAALawson}

We proceed to sketch the application of the AAA--Lawson method, which is introduced in~\cite{NT20}, to the approximation of $\ee^{\ii x}$.
We show that Proposition~\ref{prop:maincaseIII} applies in this case and the generated approximant is unitary.
The AAA--Lawson method first runs the AAA method, which is summarized
in the previous subsection.
This returns a set of support nodes $y_1,\ldots,y_m\in\R$
which is fixed for the following procedure.
The algorithm then proceeds to find a near-best approximant
using the barycentric rational representation $r_b$ given in~\eqref{eq.defr2}.
This requires determining proper coefficients
$\alpha_1,\ldots,\alpha_m$ and $\beta_1,\ldots,\beta_m$,
which is done in an iterative procedure as summarized below.
In contrast to the AAA algorithm, the support nodes are included
in the set of test nodes in the AAA--Lawson algorithm.

As an initial iteration step, the algorithm introduces weights $\mu_1,\ldots,\mu_n=1$.
\begin{enumerate}
\item[(a)] The coefficients $\alpha_j$ and $\beta_j$ are
computed to satisfy~\eqref{eq.linerrorlawson0},
where $\mu_k$ in~\eqref{eq.linerrorlawsonmatform}
refers to the current weights.
This provides an approximant $r_b$ which minimizes the respective weighted
linearized error.
\item[(b)] A new set of weights is computed.
Namely, we update the weights $\mu_k$ to $\mu_k|r_b(x_k)-\ee^{\ii x_k}|$.
%where the weights and approximant on the right-hand side of this formula refer to the previously computed weights and approximant.
In addition, weights are normalized.
Then the iteration proceeds with (a) until a given number
of iteration runs are done or other conditions are satisfied.
\end{enumerate}

In any case, the coefficients of the generated approximant $r_b$
minimize a weighted linearized error as in~\eqref{eq.linerrorlawson0}
and Proposition~\ref{prop:maincaseIII} applies in this setting,
assuming further conditions given therein hold true.
This shows the approximant generated by the AAA--Lawson method
has no poles on the real axis and is unitary.

\appendix

\section*{Appendix}
\section{Properties of singular vectors of Loewner-type matrices}\label{subsec.Asvd}
In the present section we show properties of singular vectors
of matrices related to the Loewner matrix~\eqref{eq.LoewnerMatrixform}.
Namely, we first consider matrices of the type~\eqref{eq.defA} introduced below.
In Subsection~\ref{subsec.appendixcaseIII} below, we consider singular vectors
of matrices which are related to expanded Loewner matrices, namely,
matrices of the type~\eqref{eq:defexpLoewner} introduced below.
Similar to previous sections, we let $x_1,\ldots,x_n\in\R$ and $y_1,\ldots,y_m\in\R$ denote given test and support nodes, respectively,
and we define the diagonal matrices
\begin{equation}\label{eq.defSfSfap}
S_F=\diag{\ee^{\ii x_j}}\in\C^{n\times n},~~~
\text{and}~~~S_f=\diag{\ee^{\ii y_j}}\in\C^{m\times m}.
\end{equation}
Let $M\in\R^{n\times m}$ be a given real matrix,
and let
\begin{equation}\label{eq.defA}
A = S_F M - M S_f \in\C^{n \times m}.
\end{equation}
The singular value decomposition of $A$ yields the factorization
\begin{equation}\label{eq.svdA}
A V = U S,
\end{equation}
where $V\in\C^{m \times m}$ and
$U\in\C^{n\times m}$ denote orthonormal bases of right and left singular vectors,
respectively, and
$S=\diag{\sigma_1,\ldots,\sigma_m}$ denotes a diagonal matrix of singular values.
For the singular values we assume the ordering
$$
\sigma_1 \geq \ldots \geq \sigma_m \geq 0.
$$

We proceed with some auxiliary results.
Define the diagonal matrices
\begin{subequations}\label{eq.defLK}
\begin{equation}
K=\diag{K_{11},\ldots,K_{mm}}\in\C^{m\times m},~~~\text{and}~~
R=\diag{R_{11},\ldots,R_{nn}}\in\C^{n\times n},
\end{equation}
with diagonal entries
\begin{equation}\label{eq.defKentries}
\begin{aligned}
 K_{jj} =
\left\{ \begin{array}{ll}
(1-\ee^{-\ii y_j})/|1-\ee^{-\ii y_j}|,&~~~\ee^{-\ii y_j} \neq 1 ~~\text{and},\\
~ \ii,&~~~ \text{otherwise},
\end{array}\right.
\end{aligned}
\end{equation}
where $\ii$ refers to the imaginary unit,
and
\begin{equation}\label{eq.defRentries}
\begin{aligned}
R_{kk} =
\left\{ \begin{array}{ll}
(1-\ee^{-\ii x_k})/|1-\ee^{-\ii x_k}|,&~~~\ee^{-\ii x_k} \neq 1 ~~\text{and},\\
~ \ii,&~~~ \text{otherwise}.
\end{array}\right.
\end{aligned}
\end{equation}
\end{subequations}

The diagonal entries of $K$ satisfy
\begin{subequations}
\begin{equation}\label{eq:SfKisKast01pre}
\ee^{\ii y_j}  K_{jj} = \ee^{\ii y_j} \frac{(1-\ee^{-\ii y_j})}{|1-\ee^{-\ii y_j}|}
= \frac{\ee^{\ii y_j}-1}{|1-\ee^{-\ii y_j}|} = -K_{jj}^\ast,
~~~\text{for $\ee^{-\ii y_j}\neq 1$}.
\end{equation}
The identity $\ee^{\ii y_j}  K_{jj} = -K_{jj}^\ast$
directly holds true for the case $\ee^{-\ii y_j}=1$
with $K_{jj}=\ii$.
In a similar manner, diagonal entries of $R$ satisfy
\begin{equation}\label{eq:SfKisKast02pre}
\ee^{\ii x_k}  R_{kk} =
\ee^{\ii x_k} \frac{(1-\ee^{-\ii x_k})}{|1-\ee^{-\ii x_k}|}
= \frac{\ee^{\ii x_k}-1}{|1-\ee^{-\ii x_k}|} = -R_{kk}^\ast,
~~~\text{for $\ee^{-\ii x_k}\neq 1$}.
\end{equation}
\end{subequations}
and this identity also holds true for the case $\ee^{-\ii x_k}=1$.
In matrix form,~\eqref{eq:SfKisKast01pre} implies
\begin{subequations}\label{eq:SfKisKast0}
\begin{equation}\label{eq:SfKisKast01}
S_f K = -K^\ast,
\end{equation}
where $S_f$ refers to the diagonal matrix
with diagonal entries $\ee^{\ii y_j}$.
In a similar manner,~\eqref{eq:SfKisKast02pre} % <-- maybe change the line break here such that reference (xx) is not confused with label (xx) for the arxiv version ---
implies
\begin{equation}\label{eq:SfKisKast02}
S_F R = -R^\ast,
\end{equation}
\end{subequations}
where $S_F$ refers to the diagonal matrix
with diagonal entries $\ee^{\ii x_k}$.

\begin{proposition}\label{prop.SVDidentities}
Let $\zeta\in\R^m$ be given, then a vector $w=\ii K \zeta$
with $w=(w_1,\ldots,w_m)^\top\in\C^m$ satisfies
\begin{equation}\label{eq:fjwjiswjast}
f_j w_j = w_j^\ast,~~~\text{for}~~j=1,\ldots,m,
\end{equation}
where $f_j=\ee^{\ii y_j}$.
\end{proposition}
\begin{proof}
Rewriting both sides of~\eqref{eq:fjwjiswjast} in matrix-vector form,
and substituting $w$ therein yields
$$
f_j w_j = (S_f w)_j = (\ii S_f K \zeta)_j, ~~~ \text{and} ~~
w_j^\ast = (-\ii K^\ast \zeta)_j,
$$
where $S_f$ denotes the diagonal matrix with diagonal entries $f_j$.
Making use of the identity~\eqref{eq:SfKisKast01}, we conclude~\eqref{eq:fjwjiswjast}.
\end{proof}

We proceed to rewrite the matrix-product $RAK$:
Substituting~\eqref{eq.defA} for $A$ and making use of~\eqref{eq:SfKisKast0},
we observe
\begin{equation}\label{eq.XtoXi}
R(S_F M - M S_f) K
= - R^\ast M K + R M K^\ast
= 2\ii \imag(R M K^\ast)\in\ii \R^{n\times m},
\end{equation}
where $\imag(X)$ denotes the entry-wise imaginary part of a matrix $X$.
We define
\begin{equation}\label{eq.notatildeX}
\widehat{A} = -\ii R A K\in\R^{n\times m}.
\end{equation}
The matrix $\widehat{A}$ is real due to~\eqref{eq.XtoXi}.
In particular,~\eqref{eq.XtoXi} shows
\begin{equation}\label{eq.notatildeXeq}
\widehat{A} = 2 \imag(R M K^\ast).
\end{equation}
The matrices $A$ and $\widehat{A}$ are similar up to a complex phase,
and thus, these matrices share the same set of singular values.
Let the singular value decomposition of the matrix $\widehat{A}\in\R^{n\times m}$ be given by
\begin{equation}\label{eq.svdx1}
\widehat{A} \widehat{V} = \widehat{U} S,
\end{equation}
where $\widehat{V}\in\R^{m\times m}$ and $ \widehat{U} \in\R^{n\times m}$
denote orthonormal bases of right and left singular vectors, respectively,
and the diagonal matrix $S\in\R^{m\times m}$ consists of the singular values of $\widehat{A}$.
Here, $S$ in~\eqref{eq.svdx1} is the same as that in~\eqref{eq.svdA}, since $A$ and $\widehat{A}$ are similar.
For the real matrix $\widehat{A}$, sets of real right and left singular vectors are accessible and unique up to a change of signs.
The bases of left and right singular vectors are orthonormal,~i.e., the matrices $\widehat{V}$ and $\widehat{U}$ are unitary\footnote{In this paper we consider orthogonality with respect to the Hermitian inner product $\langle v, u \rangle = v^\ast u$, so that a matrix with (complex) orthonormal columns is unitary. This is to be distinguished from the more common convention where orthogonal matrices are defined with respect to the inner product $\langle v, u \rangle = v^T u$, and real orthogonal matrices are distinct from unitary matrices.},
\begin{equation}\label{eq:hatVUortho}
\widehat{V}^\ast \widehat{V} = I
~~~\text{and}~~
\widehat{U}^\ast  \widehat{U} =I.
\end{equation}

\begin{proposition}\label{prop.svdAmainprop}
The matrices
\begin{equation}\label{eq.notatildeUV}
U =  -R^\ast \widehat{U} \in\C^{n\times m},~~~\text{and}~~~ V = \ii K\widehat{V} \in\C^{m\times m}
\end{equation}
correspond to orthonormal bases of right and left singular vectors of $A$, respectively.
\end{proposition}

\begin{proof}%[Proof of Proposition~\ref{prop.svdAmainprop}]
The matrices $K$ and $R$ given in~\eqref{eq.defLK} satisfy
\begin{equation}\label{eq.LKinnerprod}
R^\ast R = I,~~~\text{and}~~
K^\ast K = I.
\end{equation}
As in~\eqref{eq:hatVUortho},
the matrices $\widehat{U}$ and $\widehat{V}$ are unitary,
and together with~\eqref{eq.LKinnerprod} this implies that $U$ and $V$
as given in~\eqref{eq.notatildeUV} are unitary.

With~\eqref{eq.notatildeUV}, we have
$$
\widehat{U} = -RU ,~~~\text{and}~~~ \widehat{V} = -\ii K^\ast V.
$$
Substituting these identities in~\eqref{eq.svdx1}, we arrive at
$$
\ii R^\ast \widehat{A}K^\ast V = U S.
$$
Substituting the identity~\eqref{eq.notatildeX} therein, we conclude that $V$ and $U$
given in~\eqref{eq.notatildeUV} satisfy a singular value decomposition of $A$
which completes the proof.
\end{proof}

%%%%%%%%%%%%%%%%

\subsection{Auxiliary results for the non-interpolatory case}\label{subsec.appendixcaseIII}

We recall some previously introduced notation.
In the present section, the matrix $M\in\R^{n\times m}$
refers to a given real matrix
and $S_F\in\C^{n\times n}$ denotes the diagonal matrix with
diagonal entries $\ee^{\ii x_k}$ as in~\eqref{eq.defSfSfap}.
Our main interest in the present subsection
lies in the singular value decomposition of
\begin{equation}\label{eq:defexpLoewner}
[ M| -S_F M]\in\C^{n\times 2m}
\end{equation}
and properties of its right singular vectors.

Let $R$ be the diagonal matrix given as in~\eqref{eq.defLK}.
With $R S_F = -R^\ast$ as given in~\eqref{eq:SfKisKast02}, we have
$$
R[ M\,| -S_F M] = [ RM| R^\ast M].
$$
Due to $R$ being unitary,
the matrices $[ M| -S_F M]$ and $[ RM| R^\ast M]$
share the same set of right singular vectors
and singular values.

We proceed to show results in a slightly more general setting:
For a given matrix $\Theta\in\C^{n\times m}$, we consider the matrix
\begin{equation*}%\label{eq:defBTheta}
B=[\Theta| \conj(\Theta)]\in\C^{n\times 2m},
\end{equation*}
where $\conj(\Theta)$ denotes the matrix with complex conjugate entries of~$\Theta$.
Considering the present work, the case $\Theta=RM$ with
$B = [ RM| R^\ast M]$ is the most relevant one.

For the singular value decomposition of $B$ we write
\begin{equation}\label{eq:Bsvd}
B V_B = U_B S_B,
\end{equation}
where $S_B=\diag{\sigma_1,\ldots,\sigma_{2m}}\in\R^{2m \times 2m}$ is the diagonal matrix of singular values, and $V_B\in\C^{2m\times 2m}$ and $U_B\in\C^{n\times 2m}$
are the matrices of singular vectors.
For the singular values we assume the ordering
$$
\sigma_1 \geq \ldots \geq \sigma_{2m} \geq 0.
$$
To determine properties of the entry-wise complex phases of the singular vectors
of $B$, we proceed to introduce the unitary transformation\footnote{
    The notations $I$ and $I_k$ refer to the identity matrix
    throughout the present work.
    Here, $I_k$ explicitly refers to the $k\times k$-dimensional identity matrix.
}
\begin{equation}\label{eq:defQ}
Q = \frac{1}{\sqrt{2}}
\left(
\begin{array}{cc}
I_m & \ii I_m\\
I_m & -\ii I_m
\end{array}
\right)\in\C^{2m\times 2m},~~~
\text{which satisfies}~~
Q^\ast Q = I_{2m}.
\end{equation}
Furthermore, we introduce the matrix $\widehat{B}\in\R^{n\times 2m}$
as
\begin{equation}\label{eq:defBhat}
\widehat{B} = \big[\real(\Theta) \,| -\imag(\Theta)\big]\in\R^{n\times 2m}.
\end{equation}
The matrices $B$ and $\widehat{B}$ satisfy the relation
\begin{equation}\label{eq:BhattoB}
\widehat{B} = \frac{1}{\sqrt{2}} B Q.
%\text{and}~~
%B = \sqrt{2} \widehat{B} Q^\ast.
\end{equation}
Thus, the matrices $B$ and $\widehat{B}$
share the same set of singular values up to a factor $\sqrt{2}$.
We introduce the matrix
\begin{equation}\label{eq:SBhat}
\widehat{S}_B = \frac{1}{\sqrt{2}} S_B,
\end{equation}
where $S_B$ refers to the diagonal matrix of singular values of $B$
as in~\eqref{eq:Bsvd}.
The singular value decomposition of the real matrix $\widehat{B}$ yields a factorization
\begin{equation}\label{eq:Bhatsvd}
\widehat{B} \widehat{V}_B = \widehat{U}_B \widehat{S}_B,
\end{equation}
where $\widehat{U}_B\in\R^{n\times 2m}$ and $\widehat{V}_B\in\R^{2m \times 2m}$
denote real orthonormal bases of left and right singular values, respectively,
and $\widehat{S}_B$ is given in~\eqref{eq:SBhat}.

We proceed to re-scale $\widehat{V}_B$ by $Q$
to construct a basis of right singular vectors of $B$.
To this end we introduce the matrix
\begin{equation}\label{eq:VBtilde}
\widetilde{V}_B = Q \widehat{V}_B \in\C^{2m \times 2m},
\end{equation}
with $\widehat{V}_B\in\R^{2m \times 2m}$ as in~\eqref{eq:Bhatsvd}.

\begin{proposition}\label{prop:svdMSFMgen}
Let $\widehat{S}_B$~\eqref{eq:SBhat}, $\widehat{U}_B\in\R^{n\times 2m}$
and $\widehat{V}_B\in\R^{2m \times 2m}$ satisfy~\eqref{eq:Bhatsvd},
and let $\widetilde{V}_B\in\C^{2m \times 2m}$ be given as in~\eqref{eq:VBtilde}.
Then, $\widetilde{V}_B$ and $\widehat{U}_B$ correspond to orthonormal bases of right
and left singular vectors of $B$, respectively.
Namely, according to~\eqref{eq:Bsvd} we have the factorization
\begin{equation}\label{eq:BsvdviaBhat}
B \widetilde{V}_B = \widehat{U}_B S_B.
\end{equation}
\end{proposition}
\begin{proof}
Substituting~\eqref{eq:BhattoB} for $\widehat{B}$
and~\eqref{eq:SBhat} for $\widehat{S}_B$ in the factorization~\eqref{eq:Bhatsvd},
we have
$$
 B Q \widehat{V}_B = \widehat{U}_B S_B.
$$
Substituting $\widetilde{V}_B = Q \widehat{V}_B$ therein
we arrive at~\eqref{eq:BsvdviaBhat}.
Concerning the unitarity of $\widetilde{V}_B$ and $\widehat{U}_B$,
we recall that $\widehat{V}_B$ and $\widehat{U}_B$ refer to
orthonormal bases of right and left singular values of $\widehat{B}$.
Thus, $\widehat{U}_B$ is unitary, and
with $Q$ being an unitary transformation,
we also conclude that $\widetilde{V}_B = Q \widehat{V}_B$ is unitary.
\end{proof}

The right singular vector of $B$ corresponding to the smallest singular value
is of some interest in previous sections.
In the setting of Proposition~\ref{prop:svdMSFMgen},
the vectors $\widehat{\gamma}=\widehat{V}_B e_{2m}\in\R^{2m}$
and $\widetilde{\gamma}=\widetilde{V}_B e_{2m}\in\C^{2m}$
yield singular vectors of $\widehat{B}$ and $B$, respectively,
each corresponding to the smallest singular value of the respective matrix.
With $\widetilde{V}_B = Q\widehat{V}_B$ we have $\widetilde{\gamma}=Q\widehat{\gamma}$.
Substituting~\eqref{eq:defQ} for~$Q$ and
writing $\widehat{\gamma}=[\widehat{\alpha};\widehat{\beta}]\in\R^{2m}$, we arrive at
\begin{equation}\label{eq.deftildegamx}
\widetilde{\gamma} =
\frac{1}{\sqrt{2}}
\left(
\begin{array}{cc}
I_m & \ii I_m\\
I_m & -\ii I_m
\end{array}
\right)
\left(
\begin{array}{c}
\widehat{\alpha}\\
\widehat{\beta}
\end{array}
\right)
=
\frac{1}{\sqrt{2}}
\left(\begin{array}{c}
\widehat{\alpha} + \ii\widehat{\beta}\\
\widehat{\alpha} - \ii\widehat{\beta}
\end{array}\right)\in\C^{2m}.
\end{equation}

\begin{corollary}\label{cor:svabsabisequal}
Let $B=[RM | R^\ast M]$,~i.e.,~the case $\Theta=RM$.
Let $\widetilde{\gamma} = \widetilde{V}_B e_{2m}\in\C^{2m}$ where
$\widetilde{V}_B$ satisfies Proposition~\ref{prop:svdMSFMgen}.
We also write $\widetilde{\gamma}=[\widetilde{\alpha};\widetilde{\beta}]$
with $\widetilde{\alpha}=(\widetilde{\alpha}_1,\ldots,\widetilde{\alpha}_m)^\top\in\C^m$
and $\widetilde{\beta}=(\widetilde{\beta}_1,\ldots,\widetilde{\beta}_m)^\top\in\C^m$.
Then,
\begin{itemize}
\item the vector $\widetilde{\gamma}$ attains the minimum
\begin{equation}\label{eq:MSFMgammismin}
\big\|[ M| -S_F M]\widetilde{\gamma}\big\|_2 = \min_{v\in\C^{2m},~\|v\|_2 =1 } \big\|[ M| -S_F M]v \big\|_2,
\end{equation}
and
\item its entries satisfy
\begin{equation}\label{eq:inpropabsabisequal}
\widetilde{\alpha}_j = \widetilde{\beta}_j^\ast,~~~j=1,\ldots,m.
\end{equation}
\end{itemize}
\end{corollary}
\begin{proof}
As stated previously in the present subsection, the set of right singular vectors
and singular values of the matrices $[ M | -S_F M] $ and $B=[RM |R^\ast M]$
coincide.
Following Proposition~\ref{prop:svdMSFMgen},
the vector $\widetilde{\gamma}$ corresponds a right singular vector of $B$
corresponding to the singular value $\sigma_{2m}$,
which is the smallest singular value of $B$.
This carries over to $ [ M| -S_F M] $, thus,
$\|[ M| -S_F M]\widetilde{\gamma}\|_2=\sigma_{2m}$
minimizes~\eqref{eq:MSFMgammismin}.
Following~\eqref{eq.deftildegamx}, we have
$$
\widetilde{\alpha} = (\widehat{\alpha} + \ii\widehat{\beta})/\sqrt{2},
~~~\text{and}~~
\widetilde{\beta} = (\widehat{\alpha} - \ii\widehat{\beta})/\sqrt{2},
$$
where $\widehat{\alpha},\widehat{\beta}\in\R^{m}$ correspond to $\widehat{\gamma}=\widehat{V}_B e_{2m}$.
This implies the identity~\eqref{eq:inpropabsabisequal}.
\end{proof}

%%%%%%%%%%%%%%%%%%%%%%
\section{Unitarity of~\texorpdfstring{$(m-1,m-1)$}{(m-1,m-1)} rational interpolation at~\texorpdfstring{$2m-1$}{2m-1} nodes }\label{subsec.Aauxiliaryunitary}
\begin{proof}[Proof of Proposition~\ref{prop.2mm1gencaesunitary}]
The interpolation property~\eqref{eq.rinterpolexp2mm1} for $r=p/q$
implies
$$
|p(\theta_j)| = |q(\theta_j)|,~~~j=1,\ldots,2m-1.
$$
Applying Proposition~\ref{prop.aaaboundedp2} given below,
we conclude $|p(x)|=|q(x)|$ for $x\in\R$.
As a consequence, the sets of real zeros of $p$ and $q$
coincide, and assuming that $p$ and $q$ have no common zeros,
we further conclude that $r$ has no poles on the real axis.
These properties entail~$|r(x)|=1$ for $x\in\R$ which completes the proof.
\end{proof}

The proof of Proposition~\ref{prop.2mm1gencaesunitary} requires the following auxiliary result.
\begin{proposition}\label{prop.aaaboundedp2}
Let $p$ and $q$ denote polynomials of degree $\leq m-1$,
and let $\theta_1,\ldots,\theta_{2m-1}\in\R$ be distinct points with
\begin{equation}\label{eq.pqid2}
|p(\theta_j)| = |q(\theta_j)|,~~~j=1,\ldots,2m-1.
\end{equation}
Then,
\begin{equation}\label{eq.pqidstatement}
|p(x)| = |q(x)|,~~~x\in\R.
\end{equation}
\end{proposition}
\begin{proof}
Define $\chi = |p|^2 -|q|^2$.
Due to $p$ and $q$ being polynomials of degree $\leq m-1$
and $x$ being real, the functions $|p(x)|^2$
and $|q(x)|^2$ conform to polynomials of degree $\leq2m-2$.
Thus, $\chi\colon \R \to \R$ is a polynomial of degree $\leq2m-2$.
The identity~\eqref{eq.pqid2} implies that $\chi$ has
$2m-1$ distinct zeros on the real axis,
and as a consequence, $\chi$ is the zero polynomial
which further implies~\eqref{eq.pqidstatement}.
\end{proof}

%\input{cleanv0.bbl}
%\printbibliography

%%%%%%%%%%%%%%%%%%%%%%%%%

%%%%%%%%%%%%%%%%%%%

\end{document}